\documentclass{article}
\usepackage[utf8]{inputenc}
\usepackage{amssymb,amsmath}  
\usepackage[colorlinks]{hyperref}
\usepackage{fullpage}
\usepackage{cancel}
\usepackage{graphicx}

\usepackage{soul}
\setstcolor{red}

\newtheorem{theorem}{Theorem}[section]
\newtheorem{definition}[theorem]{Definition}

\newtheorem{remark}[theorem]{Remark}

\newcommand{\Lbdak}[1] {\Lambda^{#1}}
\newcommand{\Lbdaktw}[1] {\tilde\Lambda^{#1}}

\DeclareMathOperator{\diff}{d}
\newcommand{\Lie}[1]{\mathcal{L}_{#1}}
\newcommand{\intp}[1]{\mathrm{i}_{#1}}
\DeclareMathOperator{\tstar}{\tilde{\star}}

 \newcommand{\twdelta}{\tilde\delta}

\newcommand{\pp}[2]{\frac{\partial #1}{\partial #2}} 
\newcommand{\dede}[2]{\frac{\delta \!#1}{\delta \hspace{-0.15ex}#2}}
\newcommand{\twdede}[2]{\frac{\twdelta \!#1}{\delta \hspace{-0.15ex}#2}}
\newcommand{\dd}[2]{\frac{\diff \!#1}{\diff \!#2}}

\DeclareMathOperator{\Hh}{\mathcal{H}}

\DeclareMathOperator{\Ch}{\mathcal{C}}

\DeclareMathOperator{\Ah}{\mathcal{A}}

\DeclareMathOperator{\Bh}{\mathcal{B}}

\DeclareMathOperator{\Fh}{\mathcal{F}}

\DeclareMathOperator{\Mh}{\mathcal{M}}
\DeclareMathOperator{\PVh}{\mathcal{PV}}
\DeclareMathOperator{\PEh}{\mathcal{PE}}

\DeclareMathOperator{\vv}{\mathbf{v}}
\DeclareMathOperator{\uv}{\mathbf{u}}

\DeclareMathOperator{\Fv}{\mathbf{F}}

\DeclareMathOperator{\Rv}{\mathbf{R}}
\DeclareMathOperator{\hatnv}{\mathbf{\hat{n}_{s}}}
\DeclareMathOperator{\hattv}{\mathbf{\hat{t}_{s}}}
\DeclareMathOperator{\hatmv}{\mathbf{\hat{n}_{t}}}
\DeclareMathOperator{\hatsv}{\mathbf{\hat{t}_{t}}}
\DeclareMathOperator{\hatkv}{\mathbf{\hat{k}}}

\newcommand{\stf}[2]{{#1}^{#2}}
\newcommand{\twf}[2]{\tilde{#1}^{#2}}
\newcommand{\Dop}[1]{\mathbf{D}_{#1}}
\newcommand{\Dbarop}[1]{\mathbf{\bar{D}}_{#1}}
\newcommand{\Hodgeop}[1]{\mathbf{H}_{#1}}
\newcommand{\Hodgebarop}[1]{\mathbf{\bar{H}}_{#1}}
\newcommand{\Wedgeop}[2]{{#1} \boldsymbol{\wedge}_h {#2}}
\newcommand{\innerprod}[2]{\left< {#1} , {#2} \right>}
\newcommand{\topopair}[2]{\left<\left< {#1} , {#2} \right>\right>}

\newcommand{\Dopadj}[1]{\mathbf{D}_{#1}^\star}
\newcommand{\Dbaropadj}[1]{\mathbf{\bar{D}}_{#1}^\star}
\newcommand{\Hodgeopadj}[1]{\mathbf{H}_{#1}^\star}
\newcommand{\Hodgebaropadj}[1]{\mathbf{\bar{H}}_{#1}^\star}
\newcommand{\Wedgeopadj}[2]{{#1} \boldsymbol{\wedge}_h^\star {#2}}

\DeclareMathOperator{\Rop}{\mathbf{R}}
\DeclareMathOperator{\Wop}{\mathbf{W}}

\DeclareMathOperator{\Qop}{\mathbf{Q}}
\DeclareMathOperator{\Hop}{\mathbf{H}}
\DeclareMathOperator{\Iop}{\mathbf{I}}
\DeclareMathOperator{\Jop}{\mathbf{J}}
\DeclareMathOperator{\Phiop}{\boldsymbol{\phi}}
\DeclareMathOperator{\PhiTop}{\boldsymbol{\phi}^T}

\newcommand{\vform}{\twf{\mu}{n}}

\title{An interpretation of TRiSK-type schemes from a discrete exterior calculus perspective}
\author{Christopher Eldred\footnote{Center for Computing Research, Sandia National Laboratories, celdred@sandia.gov}, Werner Bauer\footnote{University of Surrey, Department of Mathematics, w.bauer@surrey.ac.uk}}
\date{October 2022}

\begin{document}

\maketitle

\tableofcontents

\abstract{TRiSK-type numerical schemes are widely used in both atmospheric and oceanic dynamical cores, due to their discrete analogues of important properties such as energy conservation and steady geostrophic modes. In this work, we show that these numerical methods are best understood as a discrete exterior calculus (DEC) scheme applied to a Hamiltonian formulation of the rotating shallow water equations based on split exterior calculus. This comprehensive description of the differential geometric underpinnings of TRiSK-type schemes completes the one started in \cite{Thuburn2012,Eldred2017}, and provides a new understanding of certain operators in TRiSK-type schemes as discrete wedge products and topological pairings from split exterior calculus. All known TRiSK-type schemes in the literature are shown to fit inside this general framework, by identifying the (implicit) choices made for various DEC operators by the different schemes. In doing so, unexplored choices and combinations are identified that might offer the possibility of fixing known issues with TRiSK-type schemes such as operator accuracy and Hollingsworth instability.}

\section{Introduction}

The rotating shallow water equations (RSWE) are a useful simplified model in the development process of geophysical fluid dynamics models, as they possess many of the same properties: conservation laws such as total energy and potential enstrophy, a (quasi-)balanced state and waves that mediate departures from this state (inertia-gravity and Rossby waves). However, their 2D nature permits much quicker testing and development of numerical schemes, since the expense of 3D modelling can be avoided. Additionally, the numerical methods used in a shallow-water model can often be partially or fully adopted for use in a 3D model.

A prominent family of numerical methods used for the numerical modelling of the RSWE are TRiSK-type schemes. These schemes are characterized by a number of features: Arakawa C-grid staggering\footnote{Fluid height at cells centers and normal component of velocity at cell edges.}, no spurious vorticity production, potential vorticity (PV) compatibility/steady geostrophic modes and various conservation laws (total mass, total circulation, total energy and/or potential enstrophy). These schemes form the basis for the numerics of both atmospheric (DYNAMICO \cite{Dubos2015}, LMD-Z \cite{Tort2015}, MPAS-A \cite{Skamarock2012}, Wavetrisk \cite{Kevlahan2019}, ICON-IAP \cite{Gassmann2013,Gassmann2018a}) and oceanic (MPAS-O \cite{Ringler2013}) models. TRiSK-type schemes originated near the dawn of atmospheric model development, starting with the celebrated Arakawa-Lamb scheme (AL81) \cite{Arakawa1981} that conserves both total energy and potential enstrophy. Similar efforts occurred in \cite{Sadourny1975}, which also conserves total energy but only conserves enstrophy in the non-divergent limit, instead of potential enstrophy. AL81 was extended to incorporate higher-order advection for the vorticity term in \cite{Takano1984}. These schemes are all restricted to uniform rectangular or lat-lon grids. This shortcoming was addressed in \cite{Thuburn2009,Ringler2010}, which generalized the approach of AL81 to arbitrary orthogonal grids\footnote{and gave rise to the name TRiSK, based on the initials of the authors of \cite{Ringler2010}.}. A partial understanding of TRiSK-type schemes in terms of discrete exterior calculus (DEC)\footnote{Discrete exterior calculus \cite{Hirani2003} is a type of structure-preserving numerical method for arbitrary polygonal grids, based on the double deRham complex \cite{Kreeft2011,Hiemstra2014} (implemented through a straight-twisted grid structure) and split exterior calculus.} appeared in \cite{Thuburn2012}, which provided the general framework for further work \cite{Weller2012a,Weller2012b,Weller2014,Thuburn2014} extending TRiSK-type schemes to arbitrary non-orthogonal grids and exploring some of the shortcomings of these schemes (such as operator accuracy). However, this work also did not recognize the presence of discrete versions of the wedge product and the topological pairing, which turn out to be key to a complete understanding of TRiSK-type schemes as a DEC. 

In a different direction, the Hamiltonian formulation of the continuous equations was exploited in \cite{Salmon2004} to develop a quasi-Hamiltonian discrete model that generalized both \cite{Arakawa1981} and \cite{Takano1984}. This quasi-Hamiltonian formulation was used to develop a multilayer TRiSK-type shallow water model in \cite{Stewart2016}. A direct extension of AL81 to the case of logically square non-orthogonal grids, not based on either quasi-Hamiltonian formulations or DEC, is found in \cite{Toy2017}. An important step towards a deeper understanding of TRiSK-type schemes occurred in \cite{Eldred2015,Eldred2017}, which unified the quasi-Hamiltonian and DEC viewpoints but still lacked the recognition of the wedge product and topological pairing.

This paper aims to fill that gap, and will present a complete understanding of TRiSK-type schemes as a type of DEC applied to a Hamiltonian formulation of the RSWE based on split exterior calculus, building on prior work in \cite{Thuburn2012,Eldred2015,Eldred2017}. Split exterior calculus is the extension of standard exterior calculus to take into account the orientation of the underlying manifold through the introduction of straight and twisted differential forms, where the latter change sign under a change of orientation. The use of split exterior calculus allows one to actually ``split'' the RSWE into topological and metric parts \cite{Bauer2016}. When viewed from a split exterior calculus perspective there is a 1-1 mapping between quantities and operators in the continuous and discrete settings, and TRiSK-type schemes are quite elegant. This is in contrast to the vector calculus perspective, where TRiSK-type schemes appear relatively complex, with many distinct operators.

Despite their popularity, TRiSK-type schemes suffer from several significant shortcomings. In particular, they have issues with operator accuracy on some grids \cite{Thuburn2014,Eldred2017,Peixoto2016,calandrini2021comparing,Turner2022}, especially in a grid refinement context and for quasi-uniform spherical grids commonly used in atmospheric and oceanic models. Additionally, TRiSK can suffer from Hollingsworth instability \cite{Hollingsworth1983,Lazic1986,Peixoto2017,Bell2017} (also known as Symmetric Instability of the Computational Kind, SICK), a non-cancellation error in the momentum equations which results in an unphysical loss of neutral stability of internal inertia-gravity waves. This instability can cause spurious energy exchanges in balanced motions \cite{Ducousso2017} and errors in the representation of deep water renewal \cite{Soontiens2017} for ocean models, along with computational instability in atmospheric models \cite{Peixoto2017,Bell2017}. A general prescription for avoiding Hollingsworth instability has yet to be found, but there are several known remedies that can alleviate it \cite{Peixoto2017,Bell2017,Gassmann2018a}. Efforts to fix the accuracy issues in TRiSK are detailed in \cite{Peixoto2016}, but the proposed fixes lose the desirable properties such as energy conservation and PV compatibility. Tantalisingly, the scheme in \cite{Toy2017} does not suffer from any accuracy issues despite being utilized on highly distorted non-orthogonal grids (but still topologically square and periodic). It is hoped that the more complete understanding of the differential geometric underpinnings of TRiSK-type schemes presented in this paper will allow the resolution of these issues in a way that does not break the desirable properties, for arbitrary grid topologies and geometries. It is anticipated that this will heavily utilize the separation of the discrete equations into topological parts (associated with the Poisson bracket) and metric parts (associated with the Hamiltonian), and is closely related to the ideas expressed in \cite{Bauer2016,Tonti2013,Tonti2014}. By categorizing existing TRiSK-type schemes in terms of the new general DEC, unexplored operator choices are identified that do not expand the stencil beyond nearest-neighbor and offer the possibility of achieving better accuracy or alleviating Hollingsworth instability; while still keeping the desirable properties. Similar efforts in the context of split finite element methods rather than DEC are detailed in \cite{Bauer2018a,BauerBehrensCotter2021}.

The rest of this paper is structured as follows: Section \ref{split-exterior} provides a brief review of split exterior calculus, focusing on new aspects compared to standard exterior calculus; Section \ref{continuous-rswe} reviews the continuous RSWE in Hamiltonian form, using both vector calculus and split exterior calculus; Section \ref{DEC} introduces a discrete exterior calculus (DEC) in 2D; Section \ref{discrete-rswe} discretizes the RSWE using this DEC; Section \ref{scheme-properties} discusses properties of the scheme and it's dependence on operator properties from DEC;  Section \ref{operator-choices} identifies specific DEC operator choices with the required properties; Section \ref{literature-schemes} explores the relationship between TRiSK-type schemes in the literature and the DEC-based scheme presented here by identifying the choices these schemes (implicitly) made for various operators in the DEC context, and includes a discussion in Section \ref{unexplored} of unexplored possibilities in choices for these operators; and finally Section \ref{conclusions} offers some conclusions. Appendix \ref{diff-properties} proves some useful properties of the discrete exterior derivative and topological pairing, Appendix \ref{leibniz-appendix} derives some relationships between various desired properties of the PV wedge product, Appendix \ref{functional-derivs-appendix} derives the functional derivatives of the discrete Hamiltonian of the RSWE, and Appendix \ref{split-vector-identities} explores the relationships between quantities and operators for split exterior calculus and vector calculus in 2D.

\section{Split Exterior Calculus for $n=2$}
\label{split-exterior}

Split exterior calculus \cite{Bauer2016} is standard exterior calculus extended with an explicit notion of orientation for the underlying manifold. It is assumed that the reader has a basic knowledge of standard exterior calculus and differential geometry, for reference see standard textbooks on differential geometry \cite{Flanders1963,Abraham2012} or the review paper at \cite{Bauer2016}. In this section only the new concepts specific to split exterior calculus are introduced: straight and twisted differential forms, the twisted Hodge star and the topological pairing.

\subsection{Manifolds, Orientations and (Twisted) Forms}
Consider a Riemannian manifold $\mathbb{M}$ with boundary $\delta \mathbb{M}$ and choice of ambient orientation $Or$, on which the space of straight $k$-forms ${\omega}^{k}$ is denoted with $\Lbdak{k}$. The opposite orientation is denoted with $-Or$. A classical example of orientation on a manifold is the choice between a right-hand or left-hand rule for the cross-product in $\mathbb{R}^3$. 

Then the space of \textit{twisted} differential $k$-forms ${\tilde\omega}^{k}$ is defined similarly to straight differential forms, except that twisted differential forms change sign under a change of orientation from $Or$ to $-Or$, while straight forms do not.\footnote{More concretely, twisted differential forms are a type of (vector)-bundle valued differential forms (BVDFs), where the bundle is the pseudoscalar bundle $\Psi$. In contrast, straight differential forms are BVDFs where the bundle is $\mathbb{R}$. Since $\Psi$ has the canonical trivializations $\Psi \otimes \Psi \rightarrow \mathbb{R}$ and $\Psi \otimes \mathbb{R} \rightarrow \Psi$, the wedge product of two twisted forms produces a straight form and the wedge product of a straight form and a twisted form produces a twisted form.} This change of sign is analogous to the situation with pseudovectors and vectors or pseudoscalars and scalars in vector calculus. 

Using both straight and twisted differential forms to describe a physical system allows a formulation that is manifestly independent of the choice of ambient orientation \cite{Bauer2016}. Since the choice of orientation is (for classical physics at least) arbitrary, this is a desirable feature. These ideas are closely related to the work of Tonti  \cite{Tonti2014,Tonti2013}, who characterized a wide range of variables for physical theories in terms of their association with oriented geometric entities: source (twisted differential forms) and configuration (straight differential forms). 

\subsection{(Twisted) Hodge Star}
Rather than the standard orientation-dependent Hodge star $\star$, an orientation independent Hodge star ($\tstar$, the twisted Hodge star) is defined through
\begin{equation}
\label{twisted-star-def}
    \stf{a}{k} \wedge \tstar \stf{b}{k} = <\stf{a}{k},\stf{b}{k}> \vform
\end{equation}
where $<>$ is the pointwise $L_2$ inner product that produces a straight $0$-form and $\vform$ is the twisted volume form, with a similar equivalent definition for $\twf{a}{k}$ and $\twf{b}{k}$. The definition (\ref{twisted-star-def}) is very similar to the one used for the standard Hodge star, but with the twisted volume form $\vform$ instead of the straight volume form. It is clear from this definition that the twisted Hodge star $\tstar$ must map between straight $k$-forms and twisted $(n-k)$-forms or twisted $k$-forms and straight $(n-k)$-forms, since $ \stf{a}{k} \wedge \tstar \stf{b}{k}$ must produce a twisted $n$-form. Additionally, it has the property $\tstar \tstar = (-1)^{k(n-k)}$, just like $\star$.

\subsection{Vector Proxies}
\label{vector-proxies}
Consider a vector $\vv$ and a pseudovector $\tilde{\vv}$. Using the  flat operator $\flat$, the interior product (contraction) $\intp{}$ and the twisted volume form $\vform$, there are four vector proxies that can be constructed from $\vv$ and $\tilde{\vv}$:
\begin{eqnarray}
\stf{v}{1} &=& \vv^\flat , \\
\twf{v}{1} &=& \tilde{\vv}^\flat , \\
\stf{v}{n-1} &=& \intp{\tilde{\vv}} \vform = \tstar \twf{v}{1}, \\
\twf{v}{n-1} &=& \intp{\vv} \twf{dV}{n} = \tstar \stf{v}{1} ,
\end{eqnarray}
where in the last two equations Hirani's formula $\intp{\mathbf{x}} \alpha = (-1)^{k(n-k)} \tstar (\tstar \alpha \wedge \mathbf{x}^\flat)$ \cite{Hirani2003} for an arbitrary $k$-form $\alpha$ has been used. From these proxies it is clear that $\tstar$ transforms (straight) twisted circulations ($1$-forms) into (twisted) straight fluxes ($(n-1)$-forms), or vice versa (sometimes with a minus sign since $\tstar \tstar = (-1)^{k(n-k)}$). A direct transformation of circulations into circulations or fluxes into fluxes (i.e. between straight and twisted $k$-forms) is accomplished through the use of the wedge product $\twf{I}{0} \wedge$, for example $\stf{x}{n-1} = \twf{I}{0} \wedge \twf{x}{n-1}$.

\begin{remark}
Note here the distinction between $1$-forms (circulations) and $(n-1)$-forms (fluxes). In 2D ($n=2$) they appear to be the same object since $n-1 = 1$. However, this obscures their true nature, and leads to convoluted notation (see for example \cite{Cotter2014}) that attempts to distinguish between them. In this work we will keep a clear distinction between the two.
\end{remark}

\subsection{Topological Pairing and Topological Functional Derivatives}
Here we introduce new concepts that extend the existing split exterior calculus as introduced in \cite{Bauer2016}, namely the 
topological pairing and a corresponding topological functional derivative, both crucial for this work. The first new concept in split exterior calculus is the \emph{topological pairing}, which is defined by:
\begin{equation}\label{toppair}
    \topopair{\stf{x}{k}}{\twf{y}{n-k}} := \int_\mathbb{M} \stf{x}{k} \wedge \twf{y}{n-k} = (-1)^{k(n-k)} \int_\mathbb{M} \twf{y}{n-k} \wedge  \stf{x}{k} = (-1)^{k(n-k)} \topopair{\twf{y}{n-k}}{\stf{x}{k}} .
\end{equation}

This pairing is purely topological, since both the wedge product and integration of a differential form do not require reference to a metric. The topological pairing can be connected to the inner product (which does require a metric, for the Hodge star $\tstar$) through:
\begin{eqnarray}
      \topopair{\stf{a}{k}}{\twf{b}{n-k}} &=& \innerprod{\stf{a}{k}}{\stf{b}{k}} , \label{inner1}\\
      \topopair{\twf{a}{k}}{\stf{b}{n-k}} &=& (-1)^{k(n-k)} \innerprod{\twf{a}{k}}{\twf{b}{k}},\label{inner2}
\end{eqnarray}
recalling the inner product of two straight $k$-forms is defined as $ \innerprod{\stf{a}{k}}{\stf{b}{k}}: = \int_\mathbb{M} \stf{a}{k} \wedge \tilde \star \stf{b}{k}$ (a definition that also holds for twisted forms).

Additionally, an integration by parts rule holds:
    \begin{eqnarray}
    \topopair{\diff \stf{x}{k-1}}{\twf{y}{n-k}} + (-1)^{k-1} \topopair{\stf{x}{k-1}}{\diff \twf{y}{n-k}} = \int_{\delta \mathbb{M}} tr (\stf{x}{k-1} \wedge \twf{y}{n-k}) ,  \\
    \topopair{\diff \twf{x}{k-1}}{\stf{y}{n-k}} + (-1)^{k-1} \topopair{\twf{x}{k-1}}{\diff \stf{y}{n-k}} = \int_{\delta \mathbb{M}} tr (\twf{x}{k-1} \wedge \stf{y}{n-k}),
    \end{eqnarray}
where $tr$ is the trace operator. In the case of a domain without boundaries (such as the ones we consider here), the boundary term on the right-hand side drops out.

Secondly, we define topological functional derivatives relative to this pairing, through
 \begin{definition}\label{def_funcder_toppair}
  The \emph{topological} functional derivatives of $\Fh[\stf{x}{k}, \twf{y}{k}]: \Lbdak{k} \times \Lbdaktw{k} \rightarrow \mathbb{R} $ with respect to the straight $k$-form $\stf{x}{k}$ and twisted $k$-form $\twf{x}{k}$ and with respect to the  topological pairing $\topopair{}{}$ are defined by
  \begin{equation}\label{equ_topFD}
  \begin{split}
   & \delta \Fh: = 
   \lim_{\epsilon \rightarrow 0} \frac{1}{\epsilon}\big( \Fh [ \stf{x}{k} + \epsilon \, \stf{\omega}{k}] - \Fh [\stf{x}{k} ]\big)
                =: \topopair{\stf{\omega}{k}}{\twdede{\Fh}{\stf{x}{k}}}
                 = \int_\mathcal{M} \stf{\omega}{k} \wedge \twdede{\Fh}{\stf{x}{k}}
                \quad \forall  \stf{\omega}{k} \in \Lbdak{k} , \\
    \text{resp.} \ &
   \delta \Fh: = 
   \lim_{\epsilon \rightarrow 0} \frac{1}{\epsilon}\big( \Fh [ \twf{y}{k} + \epsilon \, \twf{\omega}{k}] - \Fh [\twf{y}{k} ]\big)
                =: \topopair{\twf{\omega}{k}}{\twdede{\Fh}{\twf{y}{k}}}
                 = \int_\mathcal{M} \twf{\omega}{k} \wedge \twdede{\Fh}{\twf{y}{k}}
                \quad \forall  \twf{\omega}{k} \in \Lbdaktw{k}  ,
  \end{split}    
  \end{equation}
  for arbitrary test functions $\stf{\omega}{k}$ (resp. $\twf{\omega}{k}$). In particular for $\stf{x}{k} \in \Lbdak{k}$ the topological functional derivative $\twdede{\Fh}{\stf{x}{k}}  \in  \Lbdaktw{n-k}$ is a twisted $(n-k)$-form, while for $\twf{y}{k} \in \Lbdaktw{k}$ the topological functional derivative $\twdede{\Fh}{\twf{y}{k}}  \in  \Lbdak{(n-k)}$ is a straight $(n-k)$-form. 
  \end{definition}
Note that counting the tildes reveals easily the nature of the functional derivatives: an odd number of tildes give twisted differential forms, while an even number gives straight differential forms.

\section{Continuous RSWE}
\label{continuous-rswe}

We will now write the rotating shallow water equations (RSWE) in Hamiltonian form using both vector calculus and split exterior calculus. The former is mostly a review of material from \cite{Salmon2004,Eldred2015}. The two subsections are structured identically to help facilitate understanding, and enable quick comparison between these two approaches. These equations apply to an arbitrary 2D manifold without boundaries (for example, the sphere or doubly periodic plane), taking into account the fact that the operators can all be defined intrinsically and the vector fields are tangent to the manifold. The coordinate system is assumed to be an arbitrary rotating Eulerian coordinate system, with rotational velocity $\Rv$ and Coriolis parameter $f = \nabla^\perp \cdot \Rv$. This rotation is almost always chosen to correspond with the rotation of the underlying geophysical system (for example, the rotation of the Earth).

In the vector calculus form, the physical quantities are written in terms of scalars ($a$ and $b$) and vectors ($\mathbf{x}$ and $\mathbf{y}$), and the key operators are: the differential operators gradient $\nabla a$, skew-gradient $\nabla^\perp a$, divergence $\nabla \cdot \mathbf{x}$ and curl $\nabla^\perp \cdot \mathbf{x}$; the "product" operators vector transpose $\mathbf{x}^\perp$, vector dot product $\mathbf{x} \cdot \mathbf{y}$, scalar product $a b$ and mixed product $a \mathbf{y}$ or $a \mathbf{y}^\perp$; and the inner product $\innerprod{}{}$, which is defined as $\innerprod{a}{b} = \int_\Omega a b$ for scalars and $\innerprod{\mathbf{x}}{\mathbf{y}} = \int_\Omega \mathbf{x} \cdot \mathbf{y}$ for vectors. Functional derivatives are defined relative to the appropriate inner product.

\begin{remark}
For a 2D manifold embedded in $\mathbb{R}^3$ with unit normal vector $\hatkv$ (in $\mathbb{R}^3$), if we consider the vector $\mathbf{x}$ to lie in $\mathbb{R}^3$ (recalling it is still tangent to the manifold), some of these operators take simpler forms in terms of common operators (gradient $\nabla_3$ and curl $\nabla_3 \times$, along with cross product $\times$ and dot product $\cdot$) from vector calculus in $\mathbb{R}^3$: $\mathbf{x}^\perp = \hatkv \times \mathbf{x}$, $\nabla^\perp x = \hatkv \times \nabla_3 x$, $\nabla^\perp \cdot \mathbf{x} = \hatkv \cdot \nabla_3 \times \mathbf{x}$. However, all of the 2D operators can also be defined intrinsically without requiring that the manifold is embedded in $\mathbb{R}^3$ and without reference to a $\hatkv$.
\end{remark}

In the split exterior calculus form the physical quantities are written in terms of straight $k$-forms $\stf{x}{k}$ and twisted $k$-forms $\twf{x}{k}$; and the key operators are the exterior derivative $\diff$, the wedge product $\wedge$, the twisted Hodge star $\tstar$, the metric inner product $\innerprod{}{}$ and the topological pairing $\topopair{}{}$. Compared to vector calculus, the number of operators is greatly reduced. As discussed above, for split exterior calculus we use the topological functional derivative defined relative to the topological pairing, since the topological functional derivative gives a natural splitting between topological parts of the equations (the Poisson brackets) and the metric parts of the equations (the Hamiltonian and it's functional derivatives). This splitting is explored in more detail below.

See Appendix \ref{split-vector-identities} for identities relating the quantities and operators of vector calculus and split exterior calculus in 2D. In particular, scalars are either $0$-forms or $2$-forms and vectors are either $1$-forms or $(n-1)$-forms (recalling that we are in 2D and therefore $n-1 = 1$), the differential operators are composed of the exterior derivative $\diff$ plus the Hodge star $\tstar$, the product operators are composed of the wedge product $\wedge$ plus the Hodge star $\tstar$, and the inner products for scalars and vectors become either the inner product for differential forms or the topological pairing.

\subsection{Vector Calculus}
Consider the RSWE with fluid height $h$, relative velocity $\uv$, absolute velocity $\vv = \uv + \Rv$, bottom topography $h_s$, Coriolis parameter $f$, relative vorticity $\zeta = \nabla^\perp \cdot \uv$, absolute vorticity $\eta = \nabla^\perp \cdot \vv = \zeta + f$, and potential vorticity $q$ (defined through $h q = \eta$). The equations of motion for $h$ and $\vv$ are
\begin{eqnarray}
\label{veqn-vc}
\pp{\vv}{t} + f \uv^\perp + \zeta \uv^\perp + \nabla (gh + gh_s + \frac{\uv \cdot \uv}{2}) &=& 0, \\
\label{heqn-vc}
\pp{h}{t} + \nabla \cdot (h \uv) &=& 0.
\end{eqnarray}
Their linearized form (around $h=H$, $\uv = 0$) is given by
\begin{eqnarray}
\label{veqn-lin-vc}
\pp{\vv}{t} + f \uv^\perp + \nabla (gh + gh_s) &=& 0 ,\\
\label{heqn-lin-vc}
\pp{h}{t} + \nabla \cdot (H \uv) &=& 0 ,
\end{eqnarray}
where primes have been dropped for convenience. The linear equations will be useful when discretizing, to understand the behaviour of the scheme for steady geostrophic modes.

\subsubsection{Hamiltonian Form}
These equations can be put into Hamiltonian form (see \cite{Shepherd1990}) using the Hamiltonian $\Hh[\vv,h]$:
\begin{equation}
\label{hamiltonian-vc}
\mathcal{H}[\vv,h] = \frac{g}{2} \innerprod{h}{h} + g \innerprod{h}{h_s} + \innerprod{h}{\frac{\uv \cdot \uv}{2}} ,
\end{equation}
and Poisson brackets $\{ \Ah, \Bh \}$ for arbitrary functionals $\Ah, \Bh$:
\begin{equation}
\label{brackets-vc}
    \{\Ah, \Bh \} = - \left< \dede{\Ah}{\vv} ,q \dede{\Bh}{\vv}^\perp \right> - \left<  \dede{\Ah}{\vv}, \nabla \dede{\Bh}{h} \right> - \left< \dede{\Ah}{h} , \nabla \cdot \dede{\Bh}{\vv} \right> ,
\end{equation}
recalling that $q$ is defined through $hq = \nabla^\perp \cdot \vv = \eta$. The last term in $\Hh[\vv,h]$ (the kinetic energy) can also be written as $\frac{1}{2} \innerprod{h \uv}{\uv}$. The functional derivatives of $\Hh[\vv,h]$ are
\begin{eqnarray}
\label{funcderivs-vc}
\Fv := \dede{\Hh}{\vv} = h \uv  , \quad\quad\quad B := \dede{\Hh}{h} = gh + gh_s + \frac{\uv \cdot \uv}{2}  .
\end{eqnarray}
Inserting the functional derivatives (\ref{funcderivs-vc}) into the Poisson brackets (\ref{brackets-vc}) and using $\Fh = \left< \hat{x},x\right>$ for $x \in (h,\vv)$ with arbitrary test function $\hat{x}\in (\hat h,\hat \vv)$ gives the equations of motion:
\begin{eqnarray}
\label{veqn-hamil-vc}
\pp{\vv}{t} + q \Fv^\perp + \nabla B &=& 0 , \\
\label{heqn-hamil-vc}
\pp{h}{t} + \nabla \cdot \Fv &=& 0 ,
\end{eqnarray}
which upon substitution of the specific values of the functional derivatives reduce to (\ref{veqn-vc}) - (\ref{heqn-vc}).

\paragraph{Casimirs}
The general expression for the RSWE Casimirs $\Ch[\vv, h]$ is given by 
\begin{equation}
\label{casimir-vc}
    \Ch[\vv, h] = \int h f(q) ,
\end{equation}
with arbitrary function $f(q)$, which has functional derivatives
\begin{eqnarray}
\label{casimir-funcderivs-vc}
\dede{\Ch}{\vv} = - \nabla^\perp f^\prime(q)  , \quad\quad\quad \dede{\Ch}{h} = f(q) - q f^\prime(q)  ,
\end{eqnarray}
where we have used $q = \frac{\nabla^\perp \cdot \vv}{h}$ and also integrated by parts (and dropped boundary terms since we are in a domain without boundaries). Insertion of (\ref{casimir-funcderivs-vc}) into the Poisson brackets (\ref{brackets-vc}) confirms (after some algebra) that
\begin{equation}
\{ \Ch, \Ah \} = 0
\end{equation}
for any $\Ah$. Important examples of RSWE Casimirs are total mass ($f=1$), total mass-weighted potential vorticity (total circulation) ($f=q$) and total potential enstrophy ($f=\frac{q^2}{2}$).

\paragraph{Vorticity Dynamics}
The evolution equation for absolute vorticity $\eta$ is obtained by taking $\nabla^\perp \cdot$ of (\ref{veqn-hamil-vc}) yielding
\begin{equation}
\label{etaeqn-vc}
    \pp{\eta}{t} + \nabla \cdot (q \Fv) = 0 .
\end{equation}
The last term in (\ref{etaeqn-vc}) can also be written as $\nabla^\perp \cdot (q \Fv^\perp)$. Using (\ref{etaeqn-vc}) and (\ref{heqn-hamil-vc}) we can get an evolution equation for $q$ as
\begin{equation}
\label{qeqn-vc}
    \pp{q}{t} + \frac{\Fv}{h} \cdot \nabla q = 0.
\end{equation}
The last term in (\ref{qeqn-vc}) can also be written as $ \frac{\Fv^\perp}{h} \cdot \nabla^\perp q$. From (\ref{qeqn-vc}) we see that $q$ is materially conserved, since $\frac{\Fv}{h} = \uv$ is the material velocity. Vorticity dynamics will be useful for understanding the PV compatibility behaviour of the discrete scheme.

\paragraph{Linearized Equations}
Following standard procedures \cite{Shepherd1993} the Hamiltonian form can be linearized around $h=H$, $\uv = 0$ to yield for arbitrary $\Ah,\Bh$
\begin{equation}
\label{brackets-linear-vc}
    \{\Ah, \Bh \}_{lin} = - \left< \dede{\Ah}{\vv} , \frac{f}{H} \dede{\Bh}{\vv}^\perp \right> - \left<  \dede{\Ah}{\vv}, \nabla \dede{\Bh}{h} \right> - \left< \dede{\Ah}{h} , \nabla \cdot \dede{\Bh}{\vv} \right> ,
\end{equation}
and
\begin{equation}
\label{hamiltonian-linear-vc}
\mathcal{H}_{lin}[\vv,h] = \int g \frac{h^2}{2} + g h h_s + H \frac{\uv \cdot \uv}{2} .
\end{equation}
The functional derivatives are thus
\begin{eqnarray}
\label{funcderivs-lin-vc}
\Fv_{lin} := \dede{\Hh_{lin}}{\vv} = H \uv , \quad\quad\quad B_{lin} := \dede{\Hh_{lin}}{h} = gh + gh_s ,
\end{eqnarray}
with equations of motion
\begin{eqnarray}
\label{veqn-lin-hamil-vc}
\pp{\vv}{t} + \frac{f}{H} \Fv_{lin}^\perp + \nabla B_{lin} &=& 0 ,\\
\label{heqn-lin-hamil-vc}
\pp{h}{t} + \nabla \cdot \Fv_{lin} &=& 0 .
\end{eqnarray}
Insertion of (\ref{funcderivs-lin-vc}) into (\ref{veqn-lin-hamil-vc}) - (\ref{heqn-lin-hamil-vc}) gives (\ref{veqn-lin-vc}) - (\ref{heqn-lin-vc}). The linearized dynamics will be useful for understanding the linear modes behaviour of the discrete scheme.

\subsection{Split Exterior Calculus}
We will now present the same development using split exterior calculus instead of vector calculus. The two sections are structured the same, and it is highly instructive to compare them.

Consider the RSWE with fluid height twisted 2-form $\twf{h}{2}$, relative velocity straight 1-form $\stf{u}{1}$, rotational velocity straight 1-form $\stf{R}{1}$, absolute velocity straight 1-form $\stf{v}{1} = \stf{u}{1} + \stf{R}{1}$, bottom topography twisted 2-form $\twf{h_s}{2}$, Coriolis parameter straight 2-form $\stf{f}{2} = \diff \stf{R}{1}$, relative vorticity straight 2-form $\stf{\zeta}{2} = \diff \stf{u}{1}$, absolute vorticity straight 2-form $\stf{\eta}{2} = \diff \stf{v}{1} = \stf{\zeta}{2} + \stf{f}{2}$, and potential vorticity twisted 0-form $\twf{q}{0}$ (defined through $\twf{h}{2} \wedge \twf{q}{0} = \stf{\eta}{2}$). The correspondence between these objects and their vector calculus analogues is given in Table \ref{split-vec-quantities-table}. These choices of types and degrees of forms for the various physical quantities are motivated by the work of Tonti and collaborators \cite{Tonti2014,Tonti2013}, who conclusively demonstrated the association of physical quantities with oriented geometric entities, which are nothing more than straight and twisted differential forms. The same choices were made in \cite{Bauer2016,Eldred2017,Thuburn2012}. 

\begin{table}
\begin{center}
\begin{tabular}{ |c|c| } 
 \hline
 Vector Calculus & Split Exterior Calculus \\ \hline \hline
 $h$ & $\twf{h}{2}$ or $\stf{h}{0}$ \\ \hline
 $\vv$ & $\stf{v}{1}$ \\ \hline
 $\uv$ & $\twf{u}{n-1}$ or $\stf{u}{1}$ \\ \hline
 $\Rv$ & $\stf{R}{1}$ \\ \hline
 $h_s$ & $\twf{h_s}{2}$ or $\stf{h_s}{0}$ \\ \hline
 $f$   & $\stf{f}{2}$ or $\twf{f}{0}$ \\ \hline
 $\eta$ & $\stf{\eta}{2}$ or $\twf{\eta}{0}$ \\ \hline
 $\zeta$ & $\stf{\zeta}{2}$ or $\twf{\zeta}{0}$ \\ \hline
 $q$ & $\twf{q}{0}$ \\ \hline
  $\Fv$ & $\twf{F}{n-1}$ \\ \hline
 $B$ & $\stf{B}{0}$ \\ 
 \hline
\end{tabular}
\end{center}
\caption{The relationships between quantities in vector calculus and split exterior calculus forms of the RSWE. Note that many of the scalar quantities such as $h$ appear as both $0$-forms and $2$-forms. This corresponds to the fact that a standard vector calculus representation does not take into account whether a "scalar" is actually a density, and also whether it is a normal or pseudo quantity. A similar thing happens with $\uv$, which can be represented as either a circulation straight $1$-form or a flux twisted $(n-1)$-form.}
\label{split-vec-quantities-table}
\end{table}

Introduce also the auxiliary quantities fluid height straight 0-form $\stf{h}{0} = \tstar \twf{h}{2}$, bottom topography straight 0-form $\stf{h_s}{0} = \tstar \twf{h_s}{2}$, relative velocity twisted $(n\!-\!1)$-form $\twf{u}{n-1} = \tstar \stf{u}{1}$, relative vorticity twisted 0-form $\twf{\zeta}{0} = \tstar \stf{\zeta}{2}$ and Coriolis parameter twisted 0-form $\twf{f}{0} = \tstar \stf{f}{2}$, which are related to the corresponding primary quantities through the Hodge star $\tstar$.

The equations of motion for $\twf{h}{2}$ and $\stf{v}{1}$ are given by:
\begin{eqnarray}
\label{veqn-split}
\pp{\stf{v}{1}}{t} + \twf{f}{0} \wedge \twf{u}{n-1} + \twf{\zeta}{0} \wedge \twf{u}{n-1} + \diff (g \stf{h}{0} + g \stf{h_s}{0} + \tstar \frac{(\stf{u}{1} \wedge \twf{u}{n-1})}{2}) = 0 , \\
\label{heqn-split}
\pp{\twf{h}{2}}{t} + \diff (\stf{h}{0} \wedge \twf{u}{n-1}) = 0.
\end{eqnarray}
Their linearized form (around $\stf{h}{0} = H$ and $\stf{u}{1} = 0$) is
\begin{eqnarray}
\label{veqn-lin-split}
\pp{\stf{v}{1}}{t} + \twf{f}{0} \wedge \twf{u}{n-1} + \diff (g \stf{h}{0} + g \stf{h_s}{0}) = 0 , \\
\label{heqn-lin-split}
\pp{\twf{h}{2}}{t} + \diff (H \twf{u}{n-1}) = 0 , 
\end{eqnarray}
where primes have been dropped for convenience.

The equations (\ref{veqn-split}) - (\ref{heqn-lin-split}) can be obtained by applying the relationships in Appendix \ref{split-vector-identities} to (\ref{veqn-vc}) - (\ref{heqn-lin-vc}).\footnote{A more fundamental derivation \cite{Holm1998} is possible by starting with an Euler-Poincar\'e variational formulation based on semi-direct product theory and written in terms of split exterior calculus; and then applying a Legendre transform to get a Lie-Poisson Hamiltonian formulation, and finally making a change of variables from momentum to velocity to get a curl-form Hamiltonian formulation as used here.}

\subsubsection{Hamiltonian Form}
Equations~(\ref{veqn-split}) - (\ref{heqn-split}) can be put into Hamiltonian form using the Hamiltonian $\Hh[\stf{v}{1}, \twf{h}{2}]$:
\begin{equation}
\label{hamiltonian-split}
\Hh[\stf{v}{1}, \twf{h}{2}] = \frac{g}{2} \innerprod{\twf{h}{2}}{\twf{h}{2}} +  g \innerprod{\twf{h}{2}}{\twf{h_s}{2}} + \innerprod{\twf{h}{2}}{\frac{\stf{u}{1} \wedge \twf{u}{n-1}}{2}},
\end{equation}
(the ordering $\stf{u}{1} \wedge \twf{u}{n-1}$ is important due to anti-symmetry of the wedge product) and Poisson bracket for arbitrary $\Ah,\Bh$:
\begin{equation}
\label{brackets-split}
    \{\Ah, \Bh \} = - \topopair{\twdede{\Ah}{\stf{v}{1}}}{\twf{q}{0} \wedge \twdede{\Bh}{\stf{v}{1}}} - \topopair{\twdede{\Ah}{\stf{v}{1}}}{\diff \twdede{\Bh}{\twf{h}{2}}} - \topopair{\twdede{\Ah}{\twf{h}{2}}}{\diff \twdede{\Bh}{\stf{v}{1}}},
\end{equation}
where $\topopair{}{}$ is the topological pairing defined in Eqn.~\eqref{toppair}.
The last term in $\Hh[\stf{v}{1}, \twf{h}{2}]$ can also be written as $\frac{1}{2} \innerprod{\stf{u}{1}}{\stf{h}{0} \wedge \stf{u}{1}} = \frac{1}{2} \innerprod{\twf{u}{n-1}}{\stf{h}{0} \wedge \twf{u}{n-1}}$ and we recall that $\stf{v}{1} = \stf{u}{1} + \stf{R}{1}$.

A key point is that the Poisson brackets are purely topological: they involve only the wedge product, exterior derivative and topological pairing, all of which can be defined without using a metric. All of the metric information is present in the Hamiltonian, through the Hodge star that is part of the definition of the inner product.

The functional derivatives of $\Hh[\stf{v}{1}, \twf{h}{2}]$ are
\begin{equation}
\label{funcderivs-split}
\twf{F}{n-1} := \twdede{\Hh}{\stf{v}{1}} = \frac{1}{2} \left[ \stf{h}{0} \wedge \twf{u}{n-1} + \tstar (\stf{h}{0} \wedge \stf{u}{1}) \right],
\quad\quad \stf{B}{0} := \twdede{\Hh}{\twf{h}{2}} = g \stf{h}{0} + g \stf{h_s}{0} + \tstar \frac{(\stf{u}{1} \wedge \twf{u}{n-1})}{2}.
\end{equation}
In taking the first functional derivative ($\twdede{\Hh}{\stf{v}{1}}$), we have used the linearity of the Hodge star to write $\delta \twf{u}{n-1} = \delta \tstar \stf{u}{1} = \tstar \delta \stf{u}{1}$. In the continuous case, it possible to simplify $\twf{F}{n-1}$ as $\stf{h}{0} \wedge \twf{u}{n-1}$. However, this simplification will be valid only for certain choices of operators in the discrete case.

Inserting the functional derivatives (\ref{funcderivs-split}) into the Poisson brackets (\ref{brackets-split}) and using $\Fh = \left< \hat{x},x\right>$ for $x \in (\twf{h}{2},\stf{v}{1})$ with arbitrary test function $\hat{x}\in (\widehat{\twf{h}{2}},\widehat{\stf{v}{1}})$ gives the equations of motion:
\begin{eqnarray}
\label{veqn-hamil-split}
\pp{\stf{v}{1}}{t} + \twf{q}{0} \wedge \twf{F}{n-1} + \diff \stf{B}{0} = 0 ,\\
\label{heqn-hamil-split}
\pp{\twf{h}{2}}{t} + \diff \twf{F}{n-1} = 0 ,
\end{eqnarray}
which upon substitution of the specific values of the functional derivatives reduce to (\ref{veqn-split}) - (\ref{heqn-split}). 

\paragraph{Casimirs}
The general expression for RSWE Casimirs $\Ch[\stf{v}{1},\twf{h}{2}]$ is given by
\begin{equation}
\label{casimir-split}
\Ch[\stf{v}{1},\twf{h}{2}] = \innerprod{\stf{h}{0}}{F(\stf{q}{0})},
\end{equation}
where $F(\stf{q}{0})$ is an arbitrary function of potential vorticity straight $0$-form $\stf{q}{0} = \twf{I}{0} \wedge \twf{q}{0}$. The functional derivatives of $\Ch[\stf{v}{1},\twf{h}{2}]$ are
\begin{eqnarray}
\label{casimir-fd-split}
\twdede{\Ch}{\stf{v}{1}} = \diff (\twf{I}{0} \wedge F^\prime), \quad\quad
\twdede{\Ch}{\twf{h}{2}} = F - \stf{q}{0} \wedge F^\prime ,
\end{eqnarray}
where $F^\prime = \frac{dF}{d\stf{q}{0}}$. Important cases are $F=1$ (total mass), $F=\stf{q}{0}$ (total circulation) and $F=\frac{\stf{q}{0} \wedge \stf{q}{0}}{2}$ (total potential enstrophy).

\paragraph{Vorticity Dynamics}
The evolution equation for absolute vorticity $\stf{\eta}{2}$ is obtained by taking $\diff$ of (\ref{veqn-hamil-split}) yielding
\begin{equation}
\label{etaeqn-split}
    \pp{\stf{\eta}{2}}{t} = \pp{(\twf{q}{0} \wedge \twf{h}{2})}{t} = \pp{(\stf{q}{0} \wedge \stf{h}{2})}{t} = - \diff (\twf{q}{0} \wedge \twf{F}{n-1}), \\
\end{equation}
and for $\stf{h}{2} = \twf{I}{0} \wedge \twf{h}{2}$ by taking $\twf{I}{0} \wedge$ of (\ref{heqn-hamil-split}) resulting in
\begin{equation}
\label{Rheqn-split}
    \pp{\stf{h}{2}}{t} + \diff (\twf{I}{0} \wedge \twf{F}{n-1}) = 0.
\end{equation}
These can be combined with (\ref{heqn-hamil-split}) to yield evolution equations for $\twf{q}{0}$ and $\stf{q}{0}$:
\begin{equation}
\label{qtildeeqn-split}
    \pp{\twf{q}{0}}{t} + \tstar (\frac{1}{\stf{h}{0}} \wedge \twf{F}{n-1} \wedge \diff \twf{q}{0}) = 0 ,
\end{equation}
\begin{equation}
\label{qeqn-split}
    \pp{\stf{q}{0}}{t} + \tstar (\frac{1}{\stf{h}{0}} \wedge \twf{F}{n-1} \wedge \diff \stf{q}{0}) = 0 .
\end{equation}
Note that the right-hand sides of (\ref{qtildeeqn-split}) and (\ref{qeqn-split}) are just Lie derviatives, since $\Lie{\uv} \twf{q}{0} = \tstar (\frac{1}{\stf{h}{0}} \wedge \twf{F}{n-1} \wedge \diff \twf{q}{0})$, $\Lie{\uv} \stf{q}{0} = \tstar (\frac{1}{\stf{h}{0}} \wedge \twf{F}{n-1} \wedge \diff \stf{q}{0})$ and $\Lie{\uv} \stf{\eta}{2} =  \diff (\twf{I}{0} \wedge \twf{F}{n-1})$. Thus we see that $\twf{q}{0}$ and $\stf{q}{0}$ are materially conserved.

\paragraph{Linearized Equations}
Following standard procedures \cite{Shepherd1993} the Hamiltonian form can be linearized around $\stf{h}{0}=H$, $\stf{u}{1} = 0$ to yield for arbitrary $\Ah,\Bh$:
\begin{equation}
\label{brackets-lin-split}
    \{\Ah, \Bh \}_{lin} = - \topopair{\twdede{\Ah}{\stf{v}{1}}}{\frac{\twf{f}{0}}{H} \wedge \twdede{\Bh}{\stf{v}{1}}} - \topopair{\twdede{\Ah}{\stf{v}{1}}}{\diff \twdede{\Bh}{\twf{h}{2}}} - \topopair{\twdede{\Ah}{\twf{h}{2}}}{\diff \twdede{\Bh}{\stf{v}{1}}} , 
\end{equation}
and
\begin{equation}
\label{hamiltonian-lin-split}
\mathcal{H}_{lin}[\stf{v}{1},\twf{h}{2}] = \frac{g}{2} \innerprod{\twf{h}{2}}{\twf{h}{2}} +  g \innerprod{\twf{h}{2}}{ \twf{h_s}{2}} + \frac{H}{2} \innerprod{\stf{u}{1}}{\stf{u}{1}} .
\end{equation}
The functional derivatives are thus
 \begin{eqnarray}
\label{funcderivs-lin-split}
\twf{F}{n-1}_{lin} := \twdede{\Hh_{lin}}{\stf{v}{1}} = H \twf{u}{n-1} , \quad\quad\quad B_{lin} := \twdede{\Hh_{lin}}{\twf{h}{2}} = g\stf{h}{0} + g\stf{h_s}{0} ,
\end{eqnarray}
with equations of motion
\begin{eqnarray}
\label{veqn-hamil-lin-split}
\pp{\stf{v}{1}}{t} + \frac{\twf{f}{0}}{H} \wedge \twf{F}{n-1}_{lin} + \diff \stf{B}{0}_{lin} &=& 0 ,  \\
\label{heqn-hamil-lin-split}
\pp{\twf{h}{2}}{t} + \diff \twf{F}{n-1}_{lin} &=& 0 .
\end{eqnarray}
Insertion of (\ref{funcderivs-lin-split}) into (\ref{veqn-hamil-lin-split}) - (\ref{heqn-hamil-lin-split}) gives (\ref{veqn-lin-split}) - (\ref{heqn-lin-split}).

\section{A Discrete Exterior Calculus in 2D}
\label{DEC}

In order to discretize the split exterior calculus form of the RSWE, we now introduce a discrete exterior calculus (DEC) in 2D without boundaries.  Much of this material is based on \cite{Thuburn2012} and \cite{Eldred2017}, with some extensions. Fundamentally, DEC is composed of discrete versions of the key elements of split exterior calculus: straight $k$-forms $\stf{x}{k}$, twisted $k$-forms $\twf{x}{k}$, exterior derivative $\diff$, wedge product $\wedge$, Hodge star $\tstar$, inner product $\innerprod{}{}$ and the topological pairing $\topopair{}{}$ (along with functional derviatives associated to the topological pairing). As shown in Section \ref{literature-schemes}, these operators are in fact those that underlie TRiSK-type schemes. This connection was partially explored in \cite{Thuburn2012,Eldred2017}, but the importance of the wedge product and topological pairing was not recognized. In particular, what was missing was the formulation of the continuous RSWE in terms of split exterior calculus, from which the full correspondence between TRiSK-type schemes and DEC can be deduced.

In this paper we will (other than in Section \ref{operator-choices} and \ref{literature-schemes}, to facilitate comparison with the more common TRiSK notation) denote the discrete exterior derivative with $\Dop{k}$ and $\Dbarop{k}$; the discrete Hodge star with $\Hodgeop{k}$ and $\Hodgebarop{k}$ and the discrete wedge product with $\Wedgeop{}{}$ (all explained further below). The discrete inner product $\innerprod{}{}$, topological pairing $\topopair{}{}$ and topological functional derivative $\twdede{}{}$ will share the same notation as their continuous counterparts.

\subsection{Grids}
\label{grids}
For DEC, a pair of related grids are used: a straight grid and a twisted grid, each composed of $k$-cells ($k=0,1,2$). A $k$-cell is an oriented geometric entity of dimension $k$: a $0$-cell is a point ($v$ or $\tilde{v}$), a $1$-cell is an edge ($e$ or $\tilde{e}$) and a $2$-cell is a cell ($c$ or $\tilde{c}$). A collection of $k$-cells for a given $k$ is a $k$-chain. 

\begin{remark}
In the DEC literature, the terminology primal (=straight) and dual (=twisted) are usually used. However, in the TRiSK-type scheme literature, this convention is flipped and primal=twisted and dual=straight. In this paper we avoid this confusion and will consistently use straight and twisted grid.
\end{remark}

By oriented we mean that there is a value in $\{1,-1\}$ assigned to each of the $(k-1)$-cells that make up the boundary of each $k$-cell, that denotes whether this boundary $(k-1)$-cell is oriented "towards" or "away" from the $k$-cell. To each $1$-cell (edge) we associate unit normal vectors $\hatnv$ (straight) and $\hatmv$ (twisted); and unit tangent vectors $\hattv$ (straight), $\hatsv$ (twisted) such that $\hatnv \cdot \hattv = 0$ and $\hatmv \cdot \hatsv = 0$. These vectors can vary across the edge, for example in the case of a curved edge. An \textit{orthogonal} pair of grids will have $\hatnv = -\hatsv$ and $\hattv = \hatmv$, while a \textit{non-orthogonal} pair of grids will not.  An excellent description of the grid topology and orientation in DEC can be found in \cite{Tonti2014}.

From these unit vectors we can define the orientation elements $n_{ec}$, $t_{ev}$ $\tilde{n}_{\tilde{e} \tilde{c}}$ and $\tilde{t}_{\tilde{e} \tilde{v}}$, which all take values in $\{-1,1\}$. Following exterior calculus conventions, the straight grid will be \textit{inner-oriented}, and the twisted grid will be \textit{outer-oriented}. An inner-oriented grid requires only the elements of the grid to orient itself, while an outer-oriented grid requires a notion of embedding in a larger space; or of duality with another grid. In fact, given an orientation for the straight grid, the twisted grid orientation is induced by simply setting $\tilde{t}_{\tilde{e} \tilde{v}} = n_{ec}$ and $\tilde{n}_{\tilde{e} \tilde{c}} = - t_{ev}$ (where duality between straight $k$-cells and twisted $(n-k)$-cells has been used, see below). Using straight grid orientation to set twisted grid orientation ensures that the straight and twisted grid are consistently oriented.

\begin{figure}[h]
\begin{center}
\includegraphics*[scale=0.30]{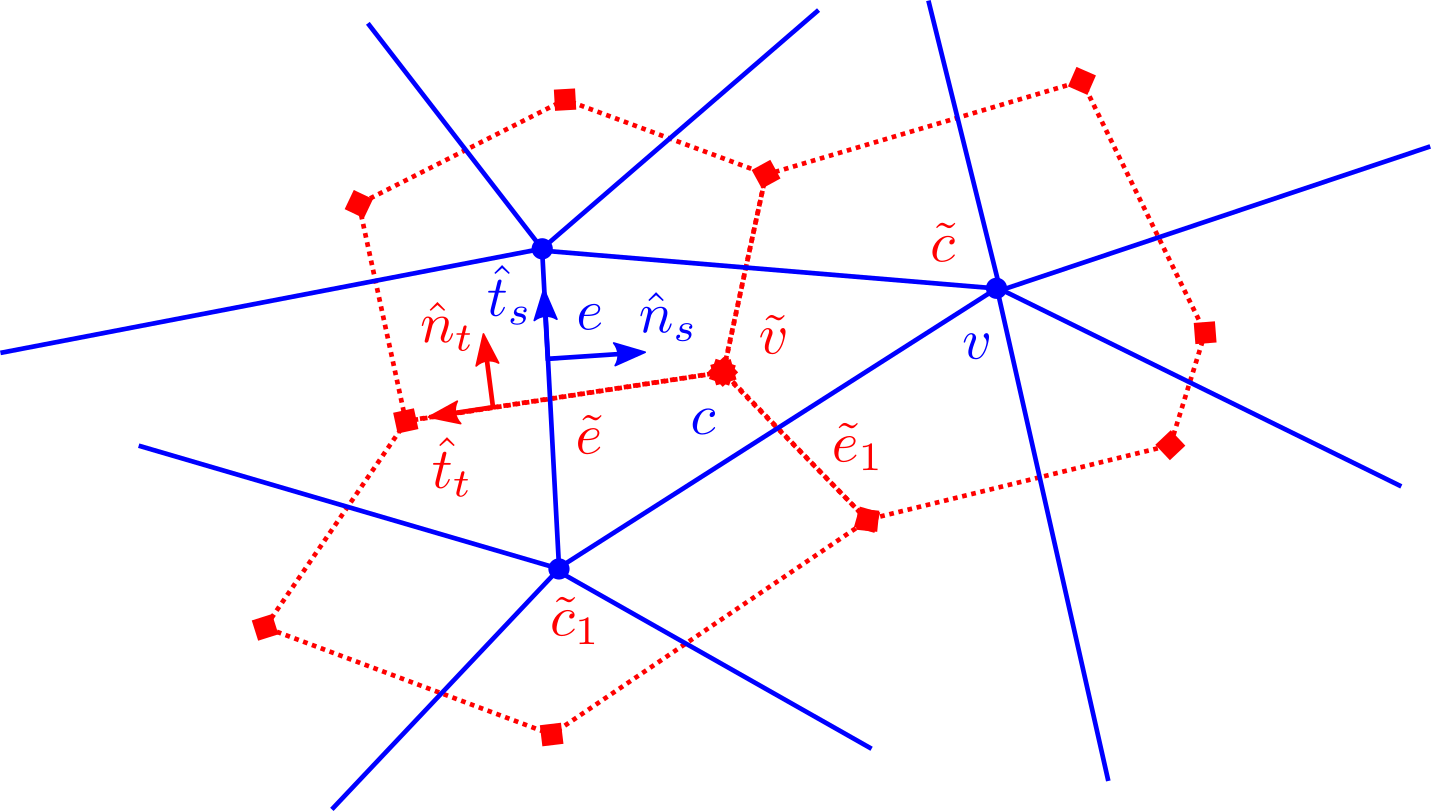}
\caption[Diagram of straight-twisted grid without boundaries]{An example straight-twisted grid without boundaries, with the straight grid in solid \textcolor{blue}{blue} and the twisted grid in dashed \textcolor{red}{red}. The unit normal vectors \textcolor{blue}{$\hatnv$} (straight) and \textcolor{red}{$\hatmv$} (twisted); and unit tangent vectors \textcolor{blue}{$\hattv$} (straight) and \textcolor{red}{$\hatsv$} (twisted) for edges $e$ and $\tilde{e}$ (which are dual to each other) are also shown. The 1-1 relationship between $k$-cells and $(n-k)$-cells is clearly visible. Some $k$-cells on both the straight and twisted grid are labeled, which correspond to the discrete differential forms in Figure \ref{grid-variables-fig}.}
\label{grid-topology-fig}
\end{center}
\end{figure}

These oriented $k$-cells are arranged into a chain complex (also known as a CW complex from algebraic topology) that can be described by a directed acyclic graph (DAG) \cite{Tonti2014}. The two grids are related through \textit{topological duality}: their respective DAG's are dual to each other in the graph theoretic sense. In particular, duality means that there is a 1-1 relationship between straight $k$-cells and twisted $(n-k)$-cells. This relationship is a fundamental part of DEC, and is used extensively (for example, in the discrete Hodge star operator). In what follows, we use $\tilde{}$ to denote the twisted grid because it will be used to represent twisted $k$-forms, and as noted before the $k$-cells themselves are denoted with $v,e,c$ (straight) or $\tilde{v},\tilde{e},\tilde{c}$ (twisted) for $k=0,1,2$, respectively. Topological duality means that there is duality between $v$ and $\tilde{c}$, $e$ and $\tilde{e}$ and $c$ and $\tilde{v}$. The topology (including orientation) for an example planar grid is shown in Figure \ref{grid-topology-fig}.

\begin{remark}
In the case that $\mathbb{M} \subset \mathbb{R}^3$, it can be useful to define a unit outward normal $\hatkv$ to $\mathbb{M}$. This can then be used to define $\hattv$ and $\hatsv$ using a right-hand rule as $\hattv = \hatkv \times \hatnv$ and $\hatsv = \hatkv \times \hatmv$. In fact, this is what is typically done for TRiSK-type schemes. This constitutes a choice of orientation for the underlying grid. However, this is not necessary and the discrete grids can be oriented without a $\hatkv$, including the case when $\mathbb{M} \not \subset \mathbb{R}^3$.
\end{remark}

On each grid, we can define a set of stencils that encompass the lower and higher-dimensional nearest-neighbor $l$-cells that surround a given $k$-cell. Specifically, for vertices $v$ we have $CV(v)$ for the nearest-neighbor cells $c$ and $EV(v)$ for the nearest-neighbor edges; for edges $e$ we have $CE(e)$ and $VE(e)$ and for cells $c$ we have $EC(c)$ and $VC(c)$. It will also be useful to introduce the stencil $ECP(e)$, which denotes the set of all edges $e^\prime$ that are in $EC(c)$ for cells $c$ in $CE(e)$; in other words it is the composition $EC(CE(e))$. Examples of these stencils can be found in \cite{Thuburn2012}.



\subsubsection{Grid Geometry}

\begin{figure}[h]
\begin{center}
\includegraphics*[scale=0.30]{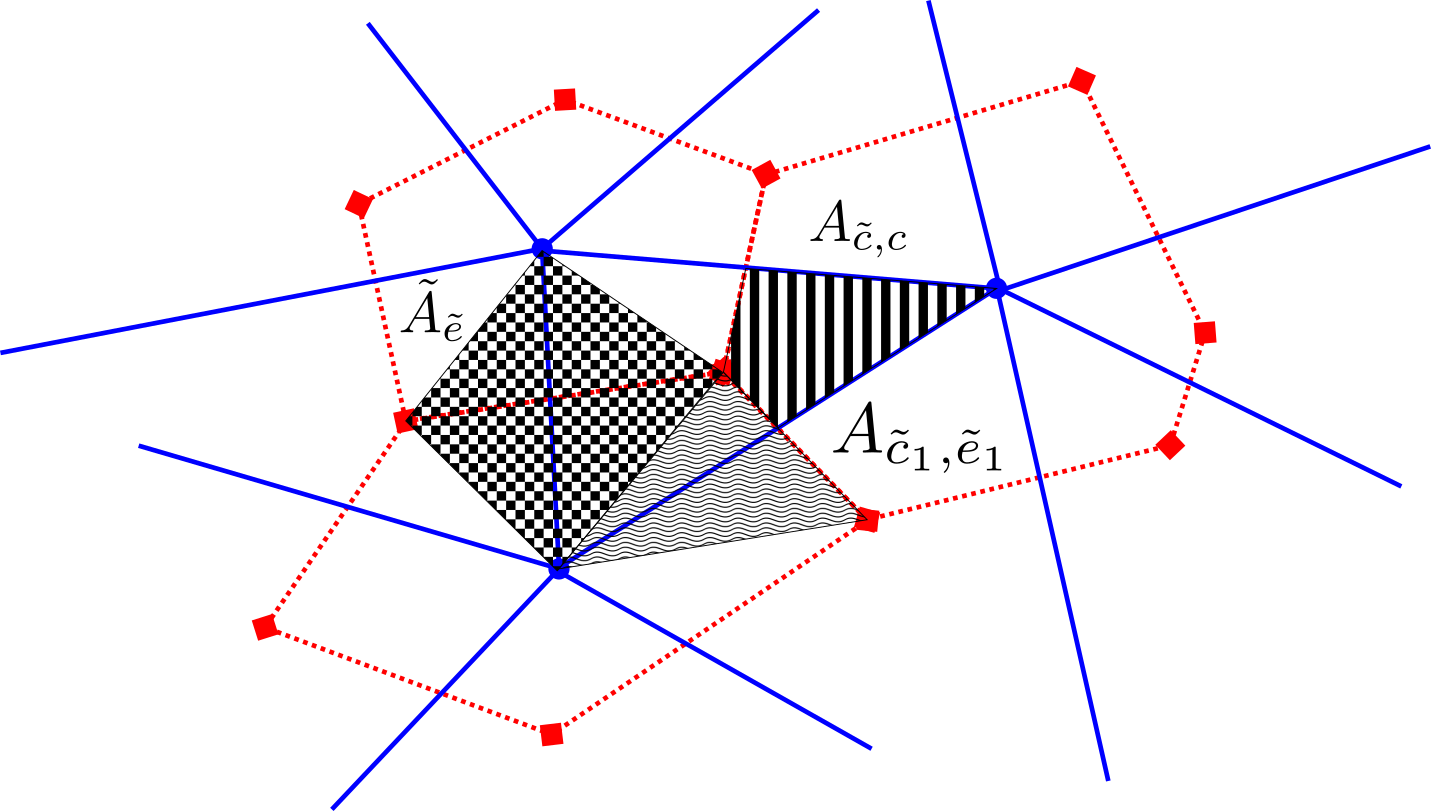}
\caption[Diagram of grid geometry]{Example geometric objects for the same straight-twisted grid as in Figure \ref{grid-topology-fig}. Specifically, the overlap areas $A_{\tilde{c},c}$ (stripes) and $A_{\tilde{c}_1,\tilde{e}_1}$ (waves) are shown, along with the extended edge area $\tilde{A}_{\tilde{e}}$ (checkerboard). These quantities are used in the definition of metric PV and KE wedge products (see definitions below). Not illustrated are the straight $k$-cells areas $A_v$ (defined to be 1), $A_e$ (length of edge $e$) and $A_c$ (area of cell $c$) and the twisted $k$-cell areas $A_{\tilde{v}}$ (defined to be 1), $A_{\tilde{e}}$ (length of edge $\tilde{e})$ and $A_{\tilde{c}}$ (area of cell $\tilde{c}$), which are self-explanatory, and are used in the Voronoi Hodge star.}
\label{grid-geom-fig}
\end{center}
\end{figure}

So far only the topological aspects of the straight and twisted grids have been discussed. To complete their description, geometric information must be introduced. Each geometric entity has a size associated with it: lengths $A_{e}$ and $A_{\tilde{e}}$ for edges, areas $A_{c}$ and $A_{\tilde{c}}$ for cells, and sizes for vertices $A_v$ and $A_{\tilde{v}}$ (defined to be $1$). It is also useful to introduce the extended area $\tilde{A}_{\tilde{e}}$ and the overlap areas $A_{c, \tilde{c}}$ and $A_{\tilde{c}, \tilde{e}}$. $\tilde{A}_{\tilde{e}}$ is the extended area associated with edge $\tilde{e}$; while $A_{c, \tilde{c}}$ is the overlap area between straight cell $c$ and twisted cell $\tilde{c}$ and $A_{\tilde{c}, \tilde{e}}$ is the overlap area between twisted cell $\tilde{c}$ and extended edge area $\tilde{A}_{\tilde{e}}$. Other overlap areas and extended areas can be similarly defined, but they are not needed for the DEC so we omit them. These sizes and overlap areas are shown in Figure \ref{grid-geom-fig}. In the same way that the topology of the twisted grid can be obtained from the topology of the straight grid through a notion of topological duality, the geometry of the twisted grid can be obtained from the straight grid through a notion of geometric duality. A commonly used approach is circumcentric/Voronoi duality, which places twisted grid vertices at the circumcenters of the straight grid cells. If the grid is optimized such that the circumcenter is also the centroid, this is a centroidal Voronoi tessellation (CVT) \cite{Lu2012,Jacobsen2013,Du2010}. Other approaches to geometric duality and grid optimization include generalized power grids \cite{engwirda2018generalised}, spring-dynamics \cite{Tomita2002,Iga2014}, barycentric \cite{Hirani2003}, Hodge-optimized triangulations \cite{mullen2011hot} and tweaking \cite{Heikes1995,Heikes2013}. The choice of geometric duality is often intimately connected to the choice of discrete Hodge star operator, for example the Voronoi Hodge star is associated with circumcentric duality and the barycentric Hodge star with barycentric duality. Certain choices of geometric duality lead to orthogonality between the straight and dual grids. It is also possible to use a grid with curved instead of straight edges: the topology does not change but the geometry does. Curved edges are useful, for example, to treat complex boundaries more accurately or to define a grid as a continuous deformation applied to a regular polygonal grid, as done in \cite{Toy2017}.
  
\subsection{Discrete Differential Forms}

\begin{figure}[h]
\begin{center}
\includegraphics*[scale=0.30]{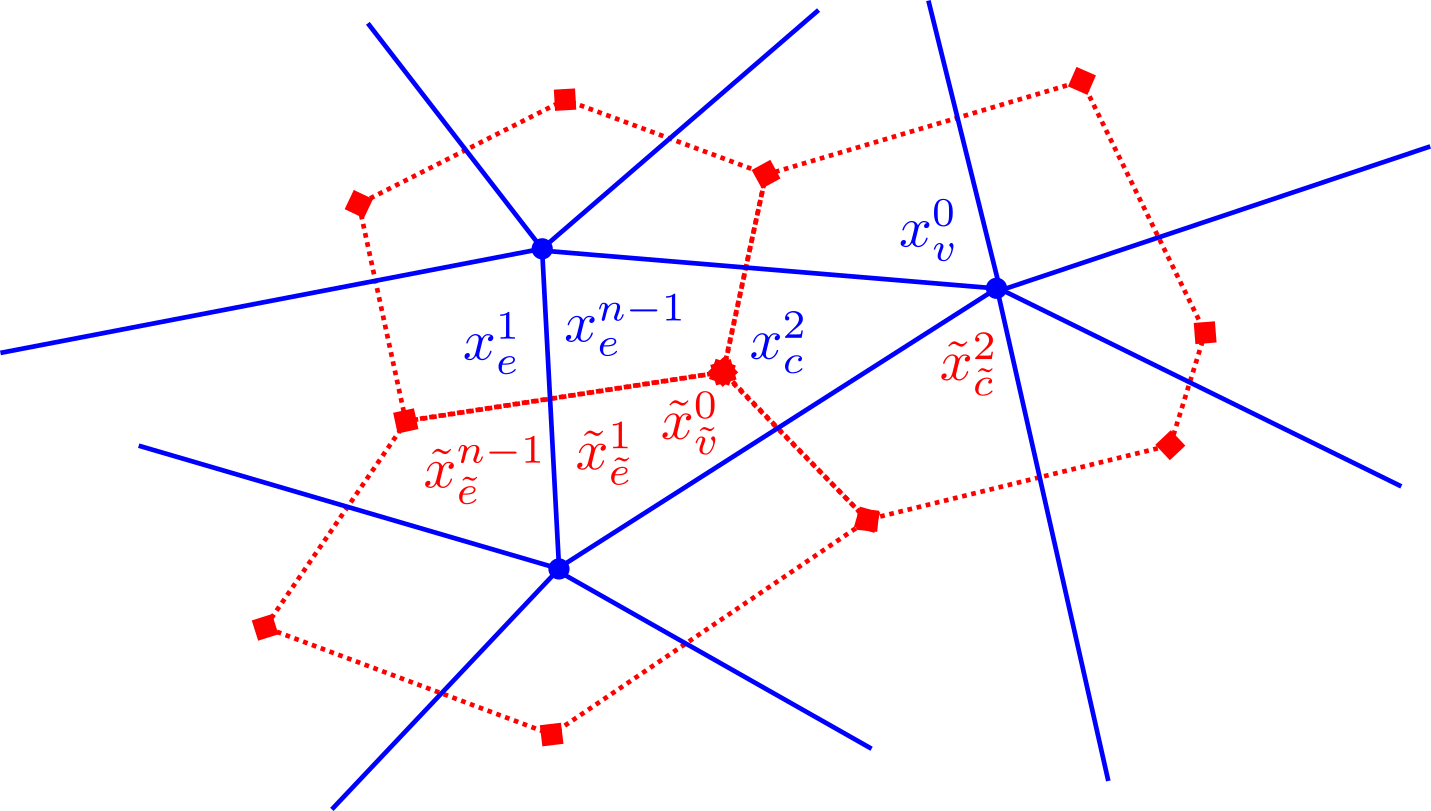}
\caption[Schematic of 2D straight-twisted grid]{Example discrete differential forms for the same straight-twisted grid as in Figure \ref{grid-topology-fig}, with the straight variables in solid \textcolor{blue}{blue} and the twisted variables in dashed \textcolor{red}{red}. Note the 1-1 duality between straight $k$-cells and twisted $(n-k)$-cells, and therefore between straight $k$-forms and twisted $(n-k)$-forms. The placement of quantities integrated over vertices ($0$-forms \textcolor{blue}{$\stf{x}{0}$} and \textcolor{red}{$\twf{x}{0}$}), edges ($1$-forms \textcolor{blue}{$\stf{x}{1}$} and \textcolor{red}{$\twf{x}{1}$} and $(n-1)$-forms \textcolor{blue}{$\stf{x}{n-1}$} and \textcolor{red}{$\twf{x}{n-1}$}) and cells ($2$-forms \textcolor{blue}{$\stf{x}{2}$} and \textcolor{red}{$\twf{x}{2}$}) is illustrated, for the same $k$-cells labeled in Figure \ref{grid-topology-fig}.}
\label{grid-variables-fig}
\end{center}
\end{figure}

A discrete differential $k$-form ($\stf{x}{k}$ or $\twf{x}{k}$)\footnote{Note that we use the same symbols for continuous and discrete quantities when there is no danger of confusion.} is represented by its action through assigning a real number to each element of a $k$-chain: this is a \textit{cochain}. In the following we will use only the notation discrete $k$-form, not $k$-cochains. We will use the straight grid for \textit{straight} (inner-oriented) forms and the twisted grid for \textit{twisted} (outer-oriented) forms. For notation, we will use $\stf{x}{k}$ (resp. $\twf{x}{k}$) to represent the column vector of discrete straight (resp. twisted) $k$-forms, and $\stf{x}{k}_p$ (resp. $\twf{x}{k}_p$) to represent the element of this vector at a specific geometric entity (where $p \in \{ v, e, c, \tilde{v}, \tilde{e}, \tilde{c}\}$ is the appropriate $k$-cell). As needed, we will denote the space of discrete straight $k$-forms with $\Lbdak{k}$ and discrete twisted $k$-forms with $\Lbdaktw{k}$, just as in the continuous case. Some example discrete $k$-forms are shown in Figure \ref{grid-variables-fig}.

More specifically, a discrete $k$-form can be thought of as the integration of the associated continuous quantity (either scalar field x or vector field $\mathbf{x}$) from vector calculus over the relevant $k$-cell:
\begin{equation}
\begin{array}{rclrcl}
\stf{x}{0}_v &=& x(v), & \twf{x}{0}_{\tilde{v}} &=&  x(\tilde{v}) , \\
\stf{x}{1}_e &=& \int_{e} (\mathbf{x} \cdot \hattv dL) , & \twf{x}{1}_{\tilde{e}} &=& \int_{\tilde{e}} (\mathbf{x} \cdot \hatsv dL) ,\\
\stf{x}{n-1}_e &=& \int_{e} (\mathbf{x} \cdot \hatnv dL), & \twf{x}{n-1}_{\tilde{e}} &=& \int_{\tilde{e}} (\mathbf{x} \cdot \hatmv dL), \\
\stf{x}{2}_c &=& \int_{c} (x dA), & \twf{x}{2}_{\tilde{c}} &=& \int_{\tilde{c}} (x dA) .
\end{array}
\end{equation}
where $dL$ and $dA$ are the differential line and area elements of edges $e$ and $\tilde{e}$ and cells $c$ and $\tilde{c}$, respectively.

As discussed in Section \ref{split-exterior},  although we are in $n=2$ and therefore $n-1=1$, using the notation $\stf{x}{1}$ and $\twf{x}{1}$ to represent both $1$-forms and $(n-1)$-forms is confusing and therefore we retain the notation $\stf{x}{n-1}$ and $\twf{x}{n-1}$ to facilitate this distinction. Using $1$ and $(n-1)$ makes clear the difference between a circulation $1$-form along an edge, and a flux $(n-1)$-form along an edge. This distinction is also seen in the unit vectors used for integration: $\hattv$ and $\hatsv$ for $1$-forms and $\hatnv$ and $\hatmv$ for $(n-1)$-forms.

\subsection{Operators}
The discrete operators are the discrete exterior derivative (denoted with $\Dop{k}$ and $\Dbarop{k}$); the discrete Hodge star (denote with $\Hodgeop{k}$ and $\Hodgebarop{k}$) and the discrete wedge product (denoted with $\Wedgeop{}{}$). For the first two we use matrix notation to indicate the representation of these operators as (sparse) matrices acting on a discrete $k$-form to produce another discrete form. In the unary operators $\Dbarop{k}$ and $\Hodgebarop{k}$ the overline indicates that the operator acts on twisted grid quantities. The wedge product $\Wedgeop{}{}$ can be represented as a (sparse) 3-tensor, since it is a binary operation that takes a discrete $k$-form and a discrete $l$-form to produce a discrete $(k+l)$-form. 

Note that the scheme presented here is very general, and to close it specific choices for the sparse coefficients and stencils appearing in the discrete Hodge stars  $\Hodgeop{k}$ and $\Hodgebarop{k}$ and wedge products $\Wedgeop{}{}$ must be made, ideally such that the properties in Section \ref{dec-properties} are obtained. These choices are usually specific to a given grid topology and geometry. Some possible choices with the desired properties are discussed in Section \ref{operator-choices}, and the specific (implicit) choices that were made by TRiSK-type schemes in the literature are identified in Section \ref{literature-schemes}.

\subsubsection{Exterior derivative}
 We will use the notation
\begin{equation}
    \Dop{k}: \Lbdak{k-1} \rightarrow \Lbdak{k} , \quad\quad \Dbarop{k}: \Lbdaktw{k-1} \rightarrow \Lbdaktw{k} ,
\end{equation}
for the discrete exterior derivative. $\Dop{k}$ and $\Dbarop{k}$ can be represented as (sparse) matrices acting on $(k-1)$-forms and producing a $k$-forms. These forms always belong to the same grid. Note that this operator exists only for $k >0$. More specifically, the discrete exterior derivatives are defined using the \textit{co-boundary operator} \cite{Tonti2014}. This definition is a \textit{combinatorial} discretization: it depends only on the topology and orientation of the underlying cell complex. Since $\diff$ is a \textit{topological} operator (one that does not require a metric), a combinatorial definition makes sense. In terms of grid orientations $n_{ec}$, $t_{ve}$, $\tilde{n}_{\tilde{e} \tilde{c}}$ and $\tilde{t}_{\tilde{e} \tilde{v}}$, the exterior derivatives are written as
\begin{eqnarray}
(\Dop{1} \stf{x}{0})_e &=& \sum_{v \in VE(e)} t_{ve} \stf{x}{0}_v , \\
(\Dop{2} \stf{x}{1})_c &=& \sum_{e \in EC(c)} n_{ec} \stf{x}{1}_e , \\
(\Dbarop{1} \twf{x}{0})_{\tilde{e}} &=& \sum_{\tilde{v} \in VE(\tilde{e})} \tilde{t}_{\tilde{v} \tilde{e}} \twf{x}{0}_{\tilde{v}} , \\
(\Dbarop{2} \twf{x}{1})_{\tilde{c}} &=& \sum_{\tilde{e} \in EC(\tilde{c})} \tilde{n}_{\tilde{e} \tilde{c}} \twf{x}{1}_{\tilde{e}} 
\end{eqnarray}

More generally, the discrete exterior derivative for a discrete $k$-form is simply the weighted sum of nearest-neighbor discrete $(k-1)$-forms, with the weights $w \in \{-1,1\}$ given by the orientation. This definition is dimension and form degree independent.

\subsubsection{Wedge product}
 We will use the notation
\begin{eqnarray}
    \Wedgeop{\stf{x}{k}}{\stf{y}{l}}: \Lbdak{k} \times \Lbdak{l} \rightarrow \Lbdak{k+l} ,\\
    \Wedgeop{\twf{x}{k}}{\stf{y}{l}}: \Lbdaktw{k} \times \Lbdak{l} \rightarrow \Lbdaktw{k+l} ,\\
    \Wedgeop{\stf{x}{k}}{\twf{y}{l}}: \Lbdak{k} \times \Lbdaktw{l} \rightarrow \Lbdaktw{k+l} ,\\
    \Wedgeop{\twf{x}{k}}{\twf{y}{l}}: \Lbdaktw{k} \times \Lbdaktw{l} \rightarrow \Lbdak{k+l} ,
\end{eqnarray}
to represent the discrete wedge product. This can be represented as a (sparse) 3-tensor acting on $k$-forms and $l$-forms to produce a $(k+l)$-forms. These forms can belong to the same grid, or different grids, depending on the specific wedge product. Note that this operator exists only when $k + l \leq n$. In algebraic topology, this operator is known as the \textit{cup product} \cite{Kotiuga2008}. As an example, consider $\stf{z}{1} = \Wedgeop{\twf{x}{0}}{\twf{y}{1}}$. An explicit form for this wedge product is 
\begin{equation}
   \stf{z}{1}_{e} = \sum_{\tilde{v}} \sum_{\tilde{e}} W_{e \tilde{v} \tilde{e}} \twf{x}{0}_{\tilde{v}} \twf{y}{1}_{\tilde{e}} ,
\end{equation}
for some set of arbitrary (sparse) coefficients $W_{e \tilde{v} \tilde{e}}$, where the target is in the first slot, and the two sources arguments are in the second and third slots. 

For TRiSK-type schemes, there are in fact only a few wedge products that must be considered. Specifically, to compute the nonlinear PV flux term and PV itself requires, respectively:
\begin{equation}
\stf{z}{1} = \Wedgeop{\twf{x}{0}}{\twf{y}{1}} , \quad\quad
\stf{z}{2} = \Wedgeop{\twf{x}{0}}{\twf{y}{2}} .
\end{equation}
To compute the kinetic energy part of the Hamiltonian and the associated mass flux and kinetic energy functional derivatives requires
\begin{equation}
\twf{z}{2} = \Wedgeop{\stf{x}{1}}{\twf{y}{1}} , \quad\quad
\twf{z}{1} = \Wedgeop{\stf{x}{0}}{\twf{y}{1}} , \quad\quad
\stf{z}{1} = \Wedgeop{\stf{x}{0}}{\stf{y}{1}} ,
\end{equation}
These cannot be chosen independently, instead one chooses one of three possible forms for the kinetic energy part of the Hamiltonian (each of which requires one of these wedge products), and then in the process of taking functional derivatives the other two wedge products arise in terms of the adjoint (see Section \ref{adjoint-ops} and Appendix \ref{functional-derivs-appendix}).

\subsubsection{Hodge star}
\label{hodge-star}
We will use the notation
\begin{equation}
\Hodgeop{k} : \Lbdak{k} \rightarrow \Lbdaktw{n-k} , \quad\quad \Hodgebarop{k} : \Lbdaktw{k} \rightarrow \Lbdak{n-k} ,
\end{equation}
to represent the discrete Hodge star. It can be represented as a (sparse) matrix acting on a $k$-form from one grid and producing a $(n-k)$-form on the other grid. By the duality between straight (twisted) $k$-cells and twisted (straight) $(n-k)$-cells, these forms have the same number of degrees of freedom. As an example, consider $\twf{y}{n-1} = \Hodgeop{1} \stf{x}{1}$. This can be written as
\begin{equation}
    \twf{y}{n-1}_{\tilde{e}} = \sum_{e} H_{\tilde{e} e} \stf{x}{1}_e
\end{equation}
for some set of arbitrary (sparse) coefficients $H_{\tilde{e} e}$, where the target is in the first slot and the source is in the second slot. 

For TRiSK-type schemes, only three discrete Hodge stars are needed:
\begin{equation}
\twf{y}{1} = \Hodgeop{1} \stf{x}{1}, \quad\quad
\stf{y}{0} = \Hodgebarop{2} \twf{x}{2}, \quad\quad
\twf{y}{0} = \Hodgeop{2} \stf{x}{2}.
\end{equation}
The remaining Hodge star operators can be implicitly defined in terms of these operators as
\begin{eqnarray}
\Hodgeop{0} := \Hodgebarop{2}^{-1} , \quad\quad
\Hodgebarop{0} := \Hodgeop{2}^{-1} , \quad\quad
\Hodgebarop{1} : = -\Hodgeop{1}^{-1} .
\end{eqnarray}
This definition is chosen so that
\begin{equation}
    \Hodgebarop{k} \Hodgeop{n-k} =  (-1)^{k(n-k)} \mathbf{I}
\end{equation}
which is a discrete analogue of $\tstar \tstar = (-1)^{k(n-k)}$. The sign in the last definition is a slight difference from \cite{Thuburn2012}. This definition requires that $\Hodgeop{1}$, $\Hodgebarop{2}$ and $\Hodgeop{2}$ are invertible, and for computational efficiency should also be sparse. Although $\Hodgeop{0}$, $\Hodgebarop{0}$ and $\Hodgebarop{1}$ will be dense (unless $\Hodgeop{1}$, $\Hodgebarop{2}$ and $\Hodgeop{2}$ are diagonal, such as in the case of the Voronoi Hodge star), this is not an issue since they are not needed in the general DEC-based TRiSK-type scheme. We could of course have started with a different set of primary Hodge stars if the scheme required it, the only requirement is that $\Hodgebarop{k} \Hodgeop{n-k} =  (-1)^{k(n-k)} \mathbf{I}$ holds.

\subsection{Topological Pairing, Functional Derivatives and Inner Products}
\label{top-pair}
It remains to define the inner product $\innerprod{}{}$, topological pairing $\topopair{}{}$ and functional derivatives (with respect to the topological pairing) $\twdede{}{}$. These definitions are done by using the discrete Hodge stars $\Hodgeop{k}$ and $\Hodgebarop{k}$, in a way that mimics $\innerprod{\stf{a}{k}}{\stf{b}{k}} = \int \stf{a}{k} \wedge \tstar \stf{b}{k}$ and $\innerprod{\twf{a}{k}}{\twf{b}{k}} = \int \twf{a}{k} \wedge \tstar \twf{b}{k}$. For all three operators, we will use the same notation for the discrete objects as their continuous counterparts.

Start with the inner product $\innerprod{}{}$ for straight and twisted forms:
\begin{equation}
\label{discrete-inner-product}
    \innerprod{\stf{x}{k}}{\stf{y}{k}} := (\stf{x}{k})^T \Hodgeop{k} \stf{y}{k}, \quad\quad \innerprod{\twf{x}{k}}{\twf{y}{k}} := (-1)^{k(n-k)} (\twf{x}{k})^T \Hodgebarop{k} \twf{y}{k} .
\end{equation}
If $\Hodgeop{k}$ and $\Hodgebarop{k}$ are positive-definite, this inner product will be positive-definitive. Some care is required with $\Hodgebarop{1}$, as it is actually negative-definite (assuming $\Hodgeop{1}$ is positive-definite) since $- \Hodgebarop{1} = \Hodgeop{1}^{-1}$. However, the minus sign for twisted 1-forms in the above definition ensure that the resulting inner product is still positive-definite.

Additionally, if $\Hodgeop{k}$ and $\Hodgebarop{k}$ are symmetric, then
\begin{eqnarray}
\label{stf-inner-symmetric}
 \innerprod{\stf{x}{k}}{\stf{y}{k}} &=&  \innerprod{\stf{y}{k}}{\stf{x}{k}} \\
 \label{twf-inner-symmetric}
 \innerprod{\twf{x}{k}}{\twf{y}{k}} &=& \innerprod{\twf{y}{k}}{\twf{x}{k}}
\end{eqnarray}
just as in the continuous case.

Now we can define the topological pairing using these inner products (similarly to ~\eqref{inner1} - \eqref{inner2}) as
\begin{align}
\label{discrete-topo-pairing-1}
\topopair{\stf{a}{k}}{\twf{b}{n-k}} &:= \innerprod{\stf{a}{k}}{\stf{b}{k}} = (\stf{a}{k})^T \Hodgeop{k} \stf{b}{k} = (\stf{a}{k})^T \twf{b}{n-k} ,\\
\label{discrete-topo-pairing-2}
\topopair{\twf{a}{k}}{\stf{b}{n-k}} \!&:= \! (-1)^{k(n-k)} \innerprod{\twf{a}{k}}{\twf{b}{k}}  = (\twf{a}{k})^T \Hodgebarop{k} \twf{b}{k} = (\twf{a}{k})^T \Hodgebarop{k} \Hodgeop{n-k} \stf{b}{n-k} =  (-1)^{k(n-k)} (\twf{a}{k})^T \stf{b}{n-k},
\end{align}
where $\twf{b}{n-k} = \Hodgeop{k} \stf{b}{k}$ and $\twf{b}{k} = \Hodgeop{n-k} \stf{b}{n-k}$. The last definition here relied on  $\Hodgebarop{k} \Hodgeop{n-k} =  (-1)^{k(n-k)} \mathbf{I}$.

This definition relies on the 1-1 relationship between straight and twisted quantities. In particular, these definitions ensure that the discrete analogues of the continuous relationships
\begin{eqnarray}
    \topopair{\stf{a}{k}}{\twf{b}{n-k}} &\stackrel{\eqref{discrete-topo-pairing-1} \& \eqref{discrete-topo-pairing-2}}{=}& (-1)^{k(n-k)} \topopair{\twf{b}{n-k}}{\stf{a}{k}}, \\
         \topopair{\stf{a}{k}}{\twf{b}{n-k}} &\stackrel{\eqref{discrete-topo-pairing-1} \& \eqref{stf-inner-symmetric}}{=}& \topopair{\stf{b}{k}}{\twf{a}{n-k}}, \\
      \topopair{\stf{a}{k}}{\twf{b}{n-k}} &\stackrel{\eqref{discrete-topo-pairing-1}}{=}& \innerprod{\stf{a}{k}}{\stf{b}{k}} ,\\
      \topopair{\twf{a}{k}}{\stf{b}{n-k}} &\stackrel{\eqref{discrete-topo-pairing-2}}{=}& (-1)^{k(n-k)} \innerprod{\twf{a}{k}}{\twf{b}{k}},
\end{eqnarray}
hold. Additionally, the topological pairing is a topological operator, so the definition used above (which is combinatorial) also makes sense. 

Finally, using this discrete topological pairing, discrete topological functional derivatives can be defined using
  \begin{eqnarray}\label{equ_topFD_discrete}
   \delta \Fh &:=& \topopair{\stf{\omega}{k}}{\twdede{\Fh}{\stf{x}{k}}} = (\stf{\omega}{k})^T \twdede{\Fh}{\stf{x}{k}} \quad \forall \stf{\omega}{k} \in \Lbdak{k} , \\
   \delta \Fh &:=& \topopair{\twf{\omega}{k}}{\twdede{\Fh}{\twf{y}{k}}} = (-1)^{k(n-k)}  (\twf{\omega}{k})^T \twdede{\Fh}{\twf{y}{k}} \quad \forall \twf{\omega}{k} \in \Lbdak{k} .
  \end{eqnarray}
As in the continuous case, a key point is that the functional derivatives lie in the dual space (defined through the Hodge star) to the variable they are taken with respect to, i.e.
\begin{equation}
\twdede{\Ah}{\stf{x}{k}} \in \Lbdaktw{n-k} , \quad\quad \twdede{\Ah}{\twf{y}{k}} \in \Lbdak{n-k} .
\end{equation}
Appendix~\ref{functional-derivs-appendix} gives an example of how this definition is used to derive functional derivatives with respect to the topological pairing.

\subsection{Operator Adjoints}
\label{adjoint-ops}

It is also useful to consider the adjoints with respect to the topological pairing of the discrete operators defined above. The adjoints will be denoted with a superscript $\star$:  $\Dopadj{k}$, $\Dbaropadj{k}$, $\Wedgeopadj{}{}$, $\Hodgeopadj{k}$ and $\Hodgebaropadj{k}$. Note that for the wedge product, there are actually two adjoints to consider since it is a binary operator. These adjoints are topological adjoints, in that they can be defined independently of a metric since they are with respect to the topological pairing.

In the continuous setting, the topological adjoint of a unary operator $\mathbf{X}$ is defined through
\begin{equation}
    \topopair{\mathbf{X}^* \stf{x}{k}}{\stf{y}{l}} := \topopair{\stf{x}{k}}{\mathbf{X} \stf{y}{l}}. 
\end{equation}
where $k$ and $l$ are defined such that the resulting topological pairing makes sense ($\stf{y}{l}$ or $\stf{x}{k}$ might even need to become twisted forms). In the discrete setting, this is written as
\begin{equation}
\label{discrete-adjoint}
    \topopair{\stf{x}{k}}{\mathbf{X} \stf{y}{l}} =  (\stf{x}{k})^T \mathbf{X} \stf{y}{l} = (\mathbf{X}^T \stf{x}{k})^T \stf{y}{l} = \topopair{\mathbf{X}^T \stf{x}{k}}{\stf{y}{l}} = \topopair{\mathbf{X}^* \stf{x}{k}}{\stf{y}{l}}, 
\end{equation}
and thus we see that $\mathbf{X}^* = \mathbf{X}^T$, or in other words the coefficients that define the discrete adjoint are just the transpose of those that define the original operator. Given $\Hodgeop{1}$, for example, $\Hodgeop{1}^*$ is defined by 
\begin{equation}
    H_{e \tilde{e}}^* = H_{\tilde{e} e} . 
\end{equation}

For a binary operator such as the wedge product, there are in fact two adjoints depending on which argument we are taking the adjoint with respect to. For example, consider the adjoint with respect to the second argument:
\begin{equation}\label{adjoint-2}
    \topopair{\twf{x}{m}}{\Wedgeop{\stf{y}{k}}{\stf{z}{l}}} = \topopair{\stf{z}{l}}{\Wedgeopadj{\stf{y}{k}}{\twf{x}{m}}} ,
\end{equation}
for $m = n-(k+l)$. In the discrete setting this is 
\begin{equation}
     \topopair{\twf{x}{m}}{\Wedgeop{\stf{y}{k}}{\stf{z}{l}}} = \sum_{m^\prime} \twf{x}{m}_{m^\prime} \sum_{k^\prime} \sum_{l^\prime} W_{m^\prime k^\prime l^\prime } \stf{y}{k}_{k^\prime} \stf{z}{l}_{l^\prime} = \sum_{l^\prime} \stf{z}{l}_{l^\prime} \sum_{m^\prime} \sum_{k^\prime} W_{l^\prime k^\prime m^\prime }^* \stf{y}{k}_{k^\prime} \twf{x}{m}_{m^\prime} = \topopair{\stf{z}{l}}{\Wedgeopadj{\stf{y}{k}}{\twf{x}{m}}} ,
\end{equation}
and thus we see that $W_{lkm}^* = W_{mkl}$, i.e. a transposition of the relevant indices of the 3-tensor.

\subsection{Operator Properties}
\label{dec-properties}
It is important for the discrete operators to have (some of) the same properties as their continuous counterparts. These properties are required in order for the discretization of the RSWE in Section \ref{discrete-rswe} to have the desirable properties discussed in Section \ref{scheme-properties}.

\subsubsection{Exterior Derivative and Topological Pairing}
The discrete exterior derivatives $\Dop{k}$ and $\Dbarop{k}$ for $k>0$ (defined using the co-boundary operator)  and the topological pairing $\topopair{}{}$ have the following set of properties:
\begin{itemize}
    \item Annihilation:
    \begin{equation}
    \label{D-annihil}
        \Dop{k} \Dop{k-1} = 0 , \quad\quad\quad \Dbarop{k} \Dbarop{k-1} = 0 ;
    \end{equation}
    Equation (\ref{D-annihil}) is the discrete analogue of $\diff \diff = 0$.
    \item Integration by parts (IBP):
    \begin{eqnarray}
    \label{D-IBP1}
    \topopair{\Dop{k} \stf{x}{k-1}}{\twf{y}{n-k}} + (-1)^{k-1} \topopair{\stf{x}{k-1}}{\Dbarop{n-k+1} \twf{y}{n-k}} = 0, \\
    \label{D-IBP2}
    \topopair{\Dbarop{k} \twf{x}{k-1}}{\stf{y}{n-k}} + (-1)^{k-1} \topopair{\twf{x}{k-1}}{\Dop{n-k+1} \stf{y}{n-k}} = 0.
    \end{eqnarray}
    which is equivalent (for the case of $n=2$) to
    \begin{equation}
    \label{D-adjoints}
        \Dbarop{2} = -\Dop{1}^T , \quad\quad\quad \Dop{2} = \Dbarop{1}^T ;
    \end{equation}
    Equations (\ref{D-IBP1}) and (\ref{D-IBP2}) are the discrete analogue of
    \begin{eqnarray}
    \topopair{\diff \stf{x}{k-1}}{\twf{y}{n-k}} + (-1)^{k-1} \topopair{\stf{x}{k-1}}{\diff \twf{y}{n-k}} = 0, \\
    \topopair{\diff \twf{x}{k-1}}{\stf{y}{n-k}} + (-1)^{k-1} \topopair{\twf{x}{k-1}}{\diff \stf{y}{n-k}} = 0.
    \end{eqnarray}
    \item Zero exterior derivative for constants: 
    \begin{equation}
    \label{D-const}
        \Dop{1} \stf{c}{0} = 0 , \quad\quad\quad \Dbarop{1} \twf{c}{0} = 0 ,
    \end{equation}
    where $\stf{c}{0}$ and $\twf{c}{0}$ are constants. Equation (\ref{D-const}) is the discrete analogue of $\diff \stf{c}{0} = 0$ and $\diff \twf{c}{0} = 0$. Note these last two properties (zero exterior derivative for constants and integration by parts) imply:
    \begin{equation}
    \label{D-stokes1}
(\twf{I}{0})^T \Dop{2} \stf{x}{1} = 0 , \quad\quad
(\stf{I}{0})^T \Dbarop{2} \twf{x}{1} = 0 ,
\end{equation}
which is the discrete analogue of Stokes theorem $\int_\Omega \diff \stf{x}{k} = \int_{\delta \Omega} \text{tr } \stf{x}{k} = 0$.
\end{itemize}

A proof of these properties for $\Dop{k}$ and $\Dbarop{k}$ based on the coboundary operator is found in Appendix \ref{diff-properties}.

\subsubsection{Wedge Product}

The continuous wedge product has three key properties:
\begin{itemize}
    \item Leibniz rule: $\diff (\stf{x}{k} \wedge \stf{y}{l}) = \diff \stf{x}{k} \wedge \stf{y}{l} + (-1)^k \stf{x}{k} \wedge \diff \stf{y}{l}$ ,
    \item Anti-symmetry: $\stf{x}{k} \wedge \stf{y}{l} = (-1)^{kl} \stf{y}{l} \wedge \stf{x}{k}$ ,
    \item Associative: $\stf{x}{m} \wedge (\stf{y}{k} \wedge \stf{z}{l}) = (\stf{x}{m} \wedge \stf{y}{k}) \wedge \stf{z}{l}$ .
\end{itemize}
Unfortunately, it has been proved using algebraic topology that it is impossible to define a discrete wedge product (i.e. a cup product) with all three properties \cite{Kotiuga2008}. However, since only the first two (anti-symmetry and the Leibniz rule) are actually needed for TRiSK-type schemes, this is not an impediment to practical usage of DEC in this case. 

In fact, only the wedge products used to compute the nonlinear PV flux term and PV itself, namely 
\begin{equation}
\stf{z}{1} = \Wedgeop{\twf{x}{0}}{\twf{y}{1}} , \quad\quad
\stf{z}{2} = \Wedgeop{\twf{x}{0}}{\twf{y}{2}} ,
\end{equation}
respectively, are required to have these properties. Specifically, we would like
\begin{itemize}
    \item Anti-symmetry for $\Wedgeop{\twf{x}{0}}{\twf{y}{1}}$: Namely, that
    \begin{equation}
    \label{Q-antisymmetry}
        \topopair{\twf{z}{1}}{\Wedgeop{\twf{x}{0}}{\twf{y}{1}}} = - \topopair{\twf{y}{1}}{\Wedgeop{\twf{x}{0}}{\twf{z}{1}}} .
    \end{equation}
    In other words, the coefficients that define $\Wedgeop{\twf{x}{0}}{\twf{y}{1}}$ must be antisymmetric in the arguments corresponding to $\twf{z}{1}$ and $\twf{y}{1}$. (\ref{Q-antisymmetry}) ensures that the (nonlinear) PV flux term conserves energy (as shown later). The vector calculus analogue of this is that $\mathbf{x} \cdot \mathbf{x}^\perp = 0$.

    \item (Partial) Leibniz rule I for $\diff$ and $\wedge$: Namely, that
    \begin{equation}
    \label{Q-pvcompat}
        \Dop{2} (\Wedgeop{\twf{I}{0}}{\twf{y}{1}}) = \Wedgeop{\twf{I}{0}}{\Dbarop{2} \twf{y}{1}} .
    \end{equation}
    Equation (\ref{Q-pvcompat}) ensures PV compatibility and steady geostrophic modes (also shown below) The vector calculus analogue of this is that $\nabla \cdot \mathbf{x} = \nabla^\perp \cdot \mathbf{x}^\perp$. 

    \item (Partial) Leibniz rule II for $\diff$ and $\wedge$ II: Namely, that 
    \begin{equation}
        \label{Q-pens}
        \Wedgeop{\twf{x}{0}}{\Dbarop{1} \twf{x}{0}} - \frac{1}{2} \Dop{1} (\Wedgeopadj{\twf{x}{0}}{\twf{x}{0}}) = 0 .
    \end{equation}
    Equation (\ref{Q-pens}) ensures potential enstrophy conservation (shown below). The vector calculus analogue of this is that $q (\nabla^\perp q)^\perp + \nabla \frac{q^2}{2} = 0$.
\item The wedge product $\stf{z}{2} = \Wedgeop{\twf{x}{0}}{\twf{y}{2}}$ has a solvability restriction on the coefficients: they must be chosen such that given $\stf{z}{2}$ and $\twf{y}{2}$ an explicit formula for $\twf{x}{0}$ is available, otherwise a linear system must be solved and the resulting scheme will be inefficient. An implicit formula for $\twf{x}{0}$ will not affect any of it's desirable properties from Section \ref{scheme-properties}, however.
\end{itemize}
An interesting question that arises is whether a connection exists between the two partial Leibniz rules (\ref{Q-pvcompat}) and (\ref{Q-pens}), motivated by the fact that in \cite{Eldred2017} a discrete wedge product that satisfied (\ref{Q-pens}) and (\ref{Q-antisymmetry}) was found to also automatically satisfy (\ref{Q-pvcompat}). This is explored in Appendix \ref{leibniz-appendix}.

\subsubsection{Hodge Star}
The discrete Hodge stars $\Hodgebarop{k}$ and $\Hodgeop{k}$ should satisfy $\Hodgebarop{n-k} \Hodgeop{k} = (-1)^{k(n-k)} \mathbf{I}$, which is ensured by the construction above in Section~\ref{hodge-star}. Additionally, $\Hodgeop{k}$ and $\Hodgebarop{k}$ must be symmetric positive-definite for the inner product $\innerprod{}{}$ to have the appropriate properties to make it an inner product, i.e. $\innerprod{a}{b} \geq 0$ with equality only if $a=0$ or $b=0$ along with $\innerprod{a}{b} = \innerprod{b}{a}$. This is clearly seen in (\ref{discrete-inner-product}).

\section{Discretization of RSWE with DEC}
\label{discrete-rswe}
Equipped with the 2D discrete exterior calculus from Section \ref{DEC}, it is straightforward to spatially discretize the split exterior calculus formulation of the RSWE from Section \ref{split-exterior}. Exactly as in the continuous case, the predicted discrete variables are the absolute velocity $\stf{v}{1}$ and fluid height $\twf{h}{2}$. The discretization simply consists in replacing $d$ with $\Dop{k}$ or $\Dbarop{k}$, $\wedge$ with $\Wedgeop{}{}$ or $\Wedgeopadj{}{}$, $\tstar$ with $\Hodgeop{k}$ or $\Hodgebarop{k}$ and using discrete versions of the inner product $\innerprod{}{}$ and topological pairing $\topopair{}{}$. The properties of these operators will ensure that the discretization inherits most of the key properties of the continuous equations.

\subsection{Discrete Hamiltonian and Functional Derivatives}
The discrete Hamiltonian $\Hh[\stf{v}{1}, \twf{h}{2}]$ is 
\begin{equation}
\Hh[\stf{v}{1}, \twf{h}{2}] = \frac{g}{2} \innerprod{\twf{h}{2}}{\twf{h}{2}} +  g \innerprod{\twf{h}{2}}{\twf{h_s}{2}} + \innerprod{\twf{h}{2}}{\frac{\Wedgeop{\stf{u}{1}}{\twf{u}{n-1}}}{2}} ,
\end{equation}
which has functional derivatives relative to the topological pairing \eqref{toppair}
\begin{equation}
\label{discrete-func-derivs}
\twf{F}{n-1} := \twdede{\Hh}{\stf{v}{1}} = \frac{1}{2} \Wedgeopadj{\stf{h}{0}}{\Hodgeop{1} \stf{u}{1}} + \frac{1}{2} \Hodgeop{1} (\Wedgeopadj{\stf{h}{0}}{\stf{u}{1}}) , \quad\quad
\stf{B}{0} := \twdede{\Hh}{\twf{h}{2}} = \frac{1}{2} \Hodgebarop{2} (\Wedgeop{\stf{u}{1}}{\twf{u}{n-1}}) + g (\stf{h}{0} + \stf{h_s}{0}) ,
\end{equation}
with topography $\twf{h_s}{2}$, recalling $\stf{h_s}{0} = \Hodgebarop{2} \twf{h_s}{2}$, $\stf{h}{0} = \Hodgebarop{2} \twf{h}{2}$ and $\twf{u}{n-1} = \Hodgeop{1} \stf{u}{1}$. The derivation of these can be found in Appendix \ref{functional-derivs-appendix}.

Thus we see that the discrete kinetic energy is $\twf{K}{2} = \frac{1}{2} (\Wedgeop{\stf{u}{1}}{\twf{u}{n-1}})$, just as in the continuous case. The mass flux is equivalent to the continuous form $\frac{1}{2} \stf{h}{0} \wedge \twf{u}{n-1} + \frac{1}{2} \tstar (\stf{h}{0} \wedge \stf{u}{1}) = \stf{h}{0} \wedge \twf{u}{n-1}$. For certain choices of $\Hodgeop{1}$ and $\Wedgeop{}{}$ this simplification is also true in the discrete setting, i.e. $\Hodgeop{1} (\Wedgeop{\stf{h}{0}}{\stf{u}{1}})  = \Wedgeop{\stf{h}{0}}{\Hodgeop{1} \stf{u}{1}}$ and then $\twf{F}{n-1} = \Wedgeop{\stf{h}{0}}{\Hodgeop{1} \stf{u}{1}} = \Wedgeop{\stf{h}{0}}{\twf{u}{n-1}}$. This is the case for the Voronoi Hodge star and the metric or combinatorial wedge products from Section \ref{scheme-properties}, as shown in Appendix \ref{simplified-form}. However, it will not be true in general.

Just as in the continuous case, the last term in $\Hh[\stf{v}{1}, \twf{h}{2}]$ (the kinetic energy) can also be written as $\frac{1}{2} \innerprod{\stf{u}{1}}{\Wedgeop{\stf{h}{0}}{\stf{u}{1}}}$ or $\frac{1}{2} \innerprod{\twf{u}{n-1}}{\Wedgeop{\stf{h}{0}}{ \twf{u}{n-1}}}$, leading respectively to slightly different functional derivatives
\begin{equation}
    \twf{F}{n-1} := \twdede{\Hh}{\stf{v}{1}} = \frac{1}{2} \Wedgeopadj{\stf{h}{0}}{\Hodgeop{1} \stf{u}{n-1}} + \frac{1}{2} \Hodgeop{1} (\Wedgeop{\stf{h}{0}}{\stf{u}{1}}) , \quad\quad
\stf{B}{0} := \twdede{\Hh}{\twf{h}{2}} = \frac{1}{2} \Hodgebarop{2} (\Wedgeopadj{\stf{u}{1}}{\twf{u}{n-1}}) + g (\stf{h}{0} + \stf{h_s}{0}) ,
\end{equation}
and
\begin{equation}
    \twf{F}{n-1} := \twdede{\Hh}{\stf{v}{1}} = \frac{1}{2} \Wedgeop{\stf{h}{0}}{\Hodgeop{1} \stf{u}{n-1}} + \frac{1}{2} \Hodgeop{1} (\Wedgeopadj{\stf{h}{0}}{\stf{u}{1}}) , \quad\quad
\stf{B}{0} := \twdede{\Hh}{\twf{h}{2}} = \frac{1}{2} \Hodgebarop{2} (\Wedgeopadj{\stf{u}{1}}{\twf{u}{n-1}}) + g (\stf{h}{0} + \stf{h_s}{0}) .
\end{equation}
These functional derivatives have the same form as (\ref{discrete-func-derivs}), the only difference is which wedge products are adjoints ($\Wedgeopadj{}{}$) and which ones are primary ($\Wedgeop{}{}$).

\subsection{Discrete Poisson Brackets}
The discrete Poisson brackets are
\begin{equation}
\label{discrete-bracket}
    \{\Ah, \Bh \} = - \topopair{\twdede{\Ah}{\stf{v}{1}}}{\Wedgeop{\twf{q}{0}}{\twdede{\Bh}{\stf{v}{1}}}} - \topopair{\twdede{\Ah}{\stf{v}{1}}}{\Dop{1} \twdede{\Bh}{\twf{h}{2}}} - \topopair{\twdede{\Ah}{\twf{h}{2}}}{\Dbarop{2} \twdede{\Bh}{\stf{v}{1}}} ,
\end{equation}
with potential vorticity $\twf{q}{0}$ defined through 
\begin{equation}
   \Wedgeop{\twf{q}{0}}{\twf{h}{2}} = \stf{\eta}{2} = \Dop{2} \stf{v}{1} = \stf{\zeta}{2} + \stf{f}{2} = \Dop{2} \stf{u}{1} + \stf{f}{2} .
\end{equation}
The definition of $\twf{q}{0}$ highlights the solvability restriction on the coefficients for the $\Wedgeop{}{}$ used to compute $\twf{q}{0}$, since $\twf{q}{0}$ should be computable without requiring a linear solve in order for the scheme to be efficient.

\subsection{Semi-discrete Equations of Motion}
The semi-discrete equations of motion that arise from these brackets and Hamiltonian are therefore
\begin{eqnarray}
\pp{\twf{h}{2}}{t} + \Dbarop{2} \twf{F}{n-1} = 0 ,\\
\pp{\stf{v}{1}}{t} + \Wedgeop{\twf{q}{0}}{\twf{F}{n-1}} + \Dop{1} \stf{B}{0} = 0 .\label{veqn}
\end{eqnarray}
These are obtained by using $\dd{\Fh}{t} = \{\Fh, \Hh \}$ for a functional $\Fh = x_p$, since in this case $\dede{\Fh}{x_p} = 1$ and $\dede{\Fh}{x} = 0$ for all other $x$; and repeating this process for all $x_p \in (\twf{h}{2}_{\tilde{c}},\stf{v}{1}_e)$. This is the discrete analogue of letting $\Fh = \left< \hat{x},x\right>$ for $x \in (\twf{h}{2},\stf{v}{1})$ with arbitrary test function $\hat{x}\in (\widehat{\twf{h}{2}},\widehat{\stf{v}{1}})$. Standard time integration schemes can then be applied to obtain the fully discrete scheme, for example explicit Runge-Kutta schemes or an (implicit) energy-conserving Poisson integrator \cite{Cohen2011} that preserves properties such as energy conservation in the fully discrete case.

\subsection{Linearized Equations}
Similar to Eqns.\eqref{veqn-hamil-split}--\eqref{heqn-hamil-split},
we linearize the above system around $\stf{h}{0} = H$ and $\stf{u}{1} = 0$ to get
\begin{eqnarray}
\pp{\stf{v}{1}}{t} + \Wedgeop{\frac{\twf{f}{0}}{H}}{\twf{F}{n-1}_{lin}}+ \Dop{1} \stf{B}{0}_{lin} &=& 0 ,\\
\pp{\twf{h}{2}}{t} + \Dbarop{2} \twf{F}{n-1}_{lin} &=& 0 ,
\end{eqnarray}
where the discrete functional derivatives read
\begin{eqnarray}
\twf{F}{n-1}_{lin} := \twdede{\Hh_{lin}}{\stf{v}{1}} = H \twf{u}{n-1}, \quad\quad\quad B_{lin} := \twdede{\Hh_{lin}}{\twf{h}{2}} = g (\stf{h}{0} + \stf{h_s}{0}),
\end{eqnarray}
for the discrete linearized Hamiltonian
\begin{equation}
\mathcal{H}_{lin}[\stf{v}{1},\twf{h}{2}] = \frac{g}{2} \innerprod{\twf{h}{2}}{\twf{h}{2}} +  g \innerprod{\twf{h}{2}}{\twf{h_s}{2}} + \frac{H}{2} \innerprod{\stf{u}{1}}{\stf{u}{1}} .
\end{equation}
As expected, almost all the wedge products (which give the nonlinearities) have disappeared, other than the linearized PV wedge product $\Wedgeop{\frac{\twf{f}{0}}{H}}{\twf{F}{n-1}_{lin}}$.

\section{Scheme Properties}
\label{scheme-properties}
The scheme presented in Section \ref{discrete-rswe} has, by construction, several desirable properties (see \cite{Staniforth2012,EldredDubos2018} for more discussion of desirable properties for schemes used in atmospheric and oceanic models). These properties come from the operator properties discussed in Section \ref{DEC} and are detailed below. They can broadly be divided into \textit{topological} properties and \textit{metric properties} (similarly to the split finite element approach taken in \cite{BauerBehrensCotter2021}). Topological properties are those that depend only on the properties of the topological operators (the exterior derivative, topological pairing and wedge product), or in other words the Poisson bracket. The metric properties involve the Hodge star and/or inner product, i.e. the Hamiltonian.

\subsection{Topological Properties}
Topological properties include structural properties such as no spurious vorticity production and PV compatibility/steady geostrophic modes; and conservation laws (total energy, mass, circulation and potential enstrophy) which arise because the discrete Poisson bracket (\ref{discrete-bracket}) is anti-symmetric and has a subset of the continuous Casimirs. Fundamentally, all of these rely on the discrete operators having discrete analogues of the key split exterior calculus identities: integration by parts for $\diff$/$\topopair{}{}$, annihilation for $\diff$, anti-symmetry for $\wedge$ and a (partial) Leibniz rule for $\diff$/$\wedge$.

\subsubsection{Total Energy Conservation}
\label{total-energy-cons}
Consider $\{\Bh, \Ah \}$:
\begin{equation}
     \{\Bh, \Ah \} = - \topopair{\twdede{\Bh}{\stf{v}{1}}}{\Wedgeop{\twf{q}{0}}{\twdede{\Ah}{\stf{v}{1}}}} - \topopair{\twdede{\Bh}{\stf{v}{1}}}{\Dop{1} \twdede{\Ah}{\twf{h}{2}}} - \topopair{\twdede{\Bh}{\twf{h}{2}}}{\Dbarop{2} \twdede{\Ah}{\stf{v}{1}}} .
\end{equation}
Using the anti-symmetry property (\ref{Q-antisymmetry}) for $\Wedgeop{}{}$ the first term can be written as
\begin{equation}
\label{Q-anti-symm}
    \topopair{\twdede{\Ah}{\stf{v}{1}}}{\Wedgeop{\twf{q}{0}}{\twdede{\Bh}{\stf{v}{1}}}} .
\end{equation}
Using the integration by parts property (\ref{D-IBP1}) for $\Dop{1}$ and $\Dbarop{2}$ the last two terms can be written as
\begin{equation}
\label{div-grad-antisymm}
    \topopair{\twdede{\Ah}{\stf{v}{1}}}{\Dop{1} \twdede{\Bh}{\twf{h}{2}}} + \topopair{\twdede{\Ah}{\twf{h}{2}}}{\Dbarop{2} \twdede{\Bh}{\stf{v}{1}}} .
\end{equation}
Therefore, combining (\ref{Q-anti-symm}) and (\ref{div-grad-antisymm}) it is clear that the discrete Poisson bracket (\ref{discrete-bracket}) is anti-symmetric:
\begin{equation}
\{\Ah, \Bh \} = -\{\Bh, \Ah \} ,
\end{equation}
and total energy is conserved since
\begin{eqnarray}
    \dd{\Hh}{t} = \{\Hh, \Hh \} = - \{\Hh, \Hh \} = 0 .
\end{eqnarray}

\subsubsection{Total Mass Conservation}
\label{total-mass-cons}
Total mass is given by
\begin{equation}
    \Mh = \innerprod{\twf{I}{2}}{\twf{h}{2}} ,
    \end{equation}
which has functional derivatives
\begin{equation}
  \twdede{\Mh}{\stf{v}{1}} = 0 , \quad\quad\quad \twdede{\Mh}{\twf{h}{2}} = \stf{I}{0} .
\end{equation}
Inserting these into (\ref{discrete-bracket}) gives
\begin{equation}
\{\Mh, \Ah \} = 0  \quad\quad \forall \Ah ,
\end{equation}
provided that $\Dop{1} \stf{I}{0} = 0$ i.e. (\ref{D-const}), and thus
\begin{eqnarray}
    \dd{\Mh}{t} = \{\Mh, \Hh \} = 0 .
\end{eqnarray}

\subsubsection{Total Circulation Conservation}
\label{total-circ-cons}
Total mass-weighted potential vorticity (i.e. total circulation) is given by
\begin{equation}
        \PVh = \innerprod{\stf{I}{2}}{\stf{\eta}{2}} \stackrel{\eqref{discrete-topo-pairing-1}}{=} \topopair{\twf{I}{0}}{\Dop{2} \stf{v}{1}} \stackrel{IBP}{=} - \topopair{\Dbarop{1} \twf{I}{0}}{\stf{v}{1}} \stackrel{\eqref{D-const}}{=} 0 ,
\end{equation}
which has functional derivatives
\begin{equation}
 \twdede{\PVh}{\stf{v}{1}} = 0 , \quad\quad\quad \twdede{\PVh}{\twf{h}{2}} = 0 .
\end{equation}
The zero functional derivatives are just a discrete analogue of the well-known fact that the total circulation for a manifold without boundaries is zero, or that the conservation of PV occurs independently of the equations of motion (it is off-shell). Therefore
\begin{equation}
\{\PVh, \Ah \} = 0  \quad\quad \forall \Ah ,
\end{equation}
and
\begin{eqnarray}
    \dd{\PVh}{t} = \{\PVh, \Hh \} = 0 .
\end{eqnarray}
Another way to see this is to simply take $\Dop{2}$ of the $\stf{v}{1}$ equation and then use the discrete Stokes theorem (\ref{D-stokes1}):
\begin{eqnarray}
    (\twf{I}{0})^T \pp{(\Dop{2} \vv)}{t} = 0
\end{eqnarray}

\subsubsection{Total Potential Enstrophy Conservation}
\label{total-pe-cons}
Total potential enstrophy is given by
\begin{equation}
    \PEh = \frac{1}{2} \innerprod{\stf{q}{2}}{\stf{\eta}{2}} 
         = \frac{1}{2} \innerprod{\Hodgebarop{0}\twf{q}{0}}{\Wedgeop{\twf{q}{0}}{\twf{h}{2}}  } 
         \stackrel{\eqref{discrete-topo-pairing-1}}{=}
         \frac{1}{2} \topopair{\twf{q}{0}}{\Wedgeop{\twf{q}{0}}{\twf{h}{2}}  } 
         \stackrel{\eqref{Q-antisymmetry}}{=}
         - \frac{1}{2} \topopair{\twf{h}{2}}{\Wedgeop{\twf{q}{0}}{\twf{q}{0}}  } 
    ,
\end{equation}
where $\stf{q}{2} = \Hodgebarop{0} \twf{q}{0}$, which has functional derivatives
\begin{eqnarray}
\twdede{\PEh}{\stf{v}{1}} = - \Dbarop{1} \twf{q}{0}, \quad\quad \twdede{\PEh}{\twf{h}{2}} = - \frac{1}{2} (\Wedgeop{\twf{q}{0}}{\twf{q}{0}}) .
\end{eqnarray}
Inserting these into (\ref{discrete-bracket}) and using the (partial) Leibniz rule II (\ref{Q-pens}) along with annihilation (\ref{D-annihil}) gives
\begin{equation}
\{\PEh, \Ah \} = 0  \quad\quad \forall \Ah ,
\end{equation}
and therefore
\begin{eqnarray}
    \dd{\PEh}{t} = \{\PEh, \Hh \} = 0 .
\end{eqnarray}
It is rare that TRiSK-type schemes choose a nonlinear PV flux operator $\Wedgeop{\twf{q}{0}}{\twf{F}{n-1}}$ that satisfies this. Instead, they usually opt for more sophisticated forms of PV advection with better behaviour for sharp gradients, such as upwinding, WENO or CLUST.

\subsubsection{No Spurious Vorticity Production}
The evolution equation for absolute vorticity $\stf{\eta}{2}$ (obtained by $\Dop{2}$ applied to (\ref{veqn})) is:
\begin{equation}
\pp{\stf{\eta}{2}}{t} + \Dop{2} (\Wedgeop{\twf{q}{0}}{\twf{F}{1}}) + \cancel{\Dop{2} \Dop{1} \stf{B}{0}} = 0 .
\end{equation}
The last term is equal to zero since $\Dop{2} \Dop{1} = 0$, i.e. (\ref{D-annihil}), which means there is no spurious vorticity production.

\subsubsection{PV Compatibility/Steady Geostrophic Modes}
\label{pv-compat-spurious-geo}
If $\twf{q}{0} = \twf{I}{0}$ (or in fact any constant), then the $\stf{\eta}{2}$ evolution equation becomes
\begin{equation}
\label{eta-discrete}
\pp{(\Wedgeop{\twf{I}{0}}{\twf{h}{2}})}{t} + \Dop{2} (\Wedgeop{\twf{I}{0}}{\twf{F}{1}}) = 0 .
\end{equation}
Now compare this to taking $\Wedgeop{\twf{I}{0}}{\square}$ of the $\twf{h}{2}$ equation
\begin{equation}
\label{Rheight-discrete}
\pp{(\Wedgeop{\twf{I}{0}}{\twf{h}{2}})}{t} + \Wedgeop{\twf{I}{0}}{\Dbarop{2} \twf{F}{1}} = 0 .
\end{equation}
If the (partial) Leibniz rule I (\ref{Q-pvcompat}) holds, then $\Dop{2} (\Wedgeop{\twf{I}{0}}{\twf{F}{1}}) = \Wedgeop{\twf{I}{0}}{\Dbarop{2} \twf{F}{1}}$ for arbitrary $\twf{F}{1}$, and (\ref{eta-discrete}) and (\ref{Rheight-discrete}) are the same equation. In that case, we say that the scheme has PV compatibility (a uniform PV field remains uniform). In the linearized equations, the same requirements arise for the presence of steady geostrophic modes (see \cite{Thuburn2012} for more details).

\subsection{Metric Properties}
The metric properties involve a combination of the topological operators along with the metric operators. They are
\begin{itemize}
    \item Linear modes: spurious stationary modes are avoided through the choice of grid staggering in the scheme (equivalent to an Arakawa C-grid) except for the Coriolis mode \cite{LeRoux2012}. The absence or presence of spurious branches of the dispersion relationship will depend on the grid topology. For example, a grid of quadrilaterals will have no spurious branches, a grid of triangles will have spurious inertia-gravity waves and a grid of hexagons will have spurious Rossby waves \cite{Randall1994,Nickovic2002,Thuburn2008}. Commonly used quasi-uniform spherical grids such as the cubed-sphere and hexagonal-pentagonal icosahedron have similar considerations, for details see \cite{Eldred2015,Thuburn2009}. The actual numerical values for the dispersion relationship itself will be a function of the (metric dependent) Hodge star operators and the choice of $\Wedgeop{\twf{I}{0}}{\twf{x}{1}}$ (the linearized version of the operator used to compute the PV flux term).
    \item Accuracy: the accuracy of the operators are determined by the exterior derivative and Hodge stars for the differential operators, and by the wedge products and Hodge stars for the product operators. This last point is especially important, since it is usually these product operators that suffer from insufficient accuracy (see Section \ref{literature-schemes} and \cite{Thuburn2014,Eldred2017}).
    \item Hollingsworth instability \cite{Hollingsworth1983,Lazic1986,Peixoto2017,Bell2017}: Hollingsworth instability occurs due to a combination of wedge products, Hodge stars and exterior derivative operators. A general prescription for avoiding it has yet to be found, but there are several known remedies that can alleviate it in the literature, all based on modifying either the kinetic energy part of the Hamiltonian \cite{Hollingsworth1983,Gassmann2013} or the nonlinear PV flux term \cite{Gassmann2018a}. Such modifications can be done in a way that preserves the other desirable properties, except potential enstrophy conservation.
\end{itemize}


\subsection{Summary of Scheme Properties}
The key features of TRiSK-type discretizations are conservation laws, no spurious vorticity production and PV compatibility/steady geostrophic modes. They all arise due to the operator properties of the topological operators (exterior derivatives, topological pairing and wedge products) and are independent of the choice of Hodge stars and grid topology and geometry. The remaining desirable properties are a good representation of linear modes, accuracy and (avoidance of) Hollingsworth instability, which interestingly are also the areas where TRiSK-type schemes can have issues. These remaining properties are largely a function of the choice of Hodge stars, specific wedge product coefficients and the choice of grid topology and geometry. This highlights the possibility of changing these choices to improve the desirable metric properties while keeping the desirable topological properties intact. This is explored more in Section \ref{unexplored}.

\section{Specific Operator Choices}
\label{operator-choices}

To close the general DEC presented in Section \ref{DEC}, choices must be made for the Hodge star and wedge product operators. Specifically, there are three choices to be made here: a choice of Hodge star, a choice of PV wedge product, and a choice of KE wedge product. These choices will determine the supported grid topologies and geometries. In this section we identify possible choices from the TRiSK, DEC and split FE literature, along with simple extensions of existing approaches to new grid topologies or geometries.

\subsection{TRiSK Notation}

To start, we will first establish the correspondence between the notation used in \cite{Eldred2015,Eldred2017} and the one used in this paper, to enable easier comparison with schemes in the literature in Section \ref{literature-schemes}. The former notation will be used only in this section and the next one.

For the Hodge star, the correspondence is
\begin{eqnarray}
    \Iop = \Hodgebarop{2} , \quad\quad
    \Hop = \Hodgeop{1} ,\quad\quad
    \Jop = \Hodgeop{2} .
\end{eqnarray}

For the PV wedge products, the correspondence is
\begin{equation}
    \Rop \twf{x}{2} = \Wedgeop{\twf{I}{0}}{\twf{x}{2}} ,\quad\quad
    \Wop \twf{x}{1} = \Wedgeop{\twf{I}{0}}{\twf{x}{1}} ,\quad\quad
    \Qop \twf{x}{1} = \Wedgeop{\twf{q}{0}}{\twf{x}{1}} ,
\end{equation}
This form of $\Rop$ is written such that an explicit formula for $\twf{x}{0}$ always exists, which is true for all the wedge products discussed here. The conditions on the PV wedge product required for energy conservation (i.e. antisymmetry (\ref{Q-antisymmetry})) can be written as
\begin{equation}
\Wop = -\Wop^T , \quad\quad\quad \Qop = - \Qop^T ,
\end{equation}
those for PV compatibility (partial Leibnitz rule I (\ref{Q-pvcompat})) as
\begin{equation}
    \Rop \Dbarop{2} = \Dop{2} \Wop ,
\end{equation}
and those for potential enstrophy conservation (partial Leibnitz rule II  (\ref{Q-pens})) as
\begin{equation}
    \Qop \Dbarop{1} \twf{q}{0} - \Dop{1} \Rop^T \frac{\twf{q}{0}}{2} = 0 .
\end{equation}
In the last equation, $\Rop^T \frac{\twf{q}{0}}{2} = \Wedgeopadj{\twf{q}{0}}{\twf{q}{0}}$, where the adjoint is from $\Wedgeop{\twf{x}{0}}{\twf{x}{2}}$ i.e. $\Rop$. These are the same formulas given in \cite{Thuburn2012,Eldred2015,Eldred2017}, up to different sign conventions.

For the kinetic energy, we will use the approach of starting with a definition for $\twf{K}{2} = \Wedgeop{\stf{x}{1}}{\twf{y}{1}}$, and then determining the wedge products needed for the mass flux as adjoint operators. The operator $\Wedgeop{\stf{x}{1}}{\twf{y}{1}}$ is closely related to the $T$ operator in \cite{Thuburn2012}, but acts on $\stf{x}{1}$ and $\twf{y}{1}$, while $T$ acts on $\stf{x}{1}$ and $\stf{y}{1}$. We do not use the $\PhiTop$ and $\Phiop$ notation from \cite{Eldred2017}, since it is not general enough to encompass the two different mass flux wedge products that arise in the case of non-diagonal $\Hop$.

\subsection{Choices for Hodge Stars: $\Iop$, $\Jop$ and $\Hop$}
A wide range of Hodge stars exist in the literature. All of the Hodge stars discussed here are low-order (between 1st and 2nd order) accurate and geometrically inflexible: they impose some sort of restriction on the topology and/or geometry of the grids. Generally speaking, they can be categorized into three types:
\begin{itemize}
    \item Voronoi Hodge star\footnote{also known as the circumcentric or diagonal Hodge star.} \cite{Hirani2013,Hirani2018}. This Hodge star is based on the 1-1 relationship between $k$-cells and $(n-k)$-cells on opposite grids, and has (diagonal) stencil coefficients equal to the ratio of areas for these two geometric entities:
\begin{equation}
I_{\tilde{c},v} = \frac{1}{A_{\tilde{c}}} , \quad\quad
J_{c,\tilde{v}} = \frac{1}{A_{c}} , \quad\quad
H_{\tilde{e},e} = \frac{A_{e}}{A_{\tilde{e}}} .
\end{equation}
More generally, it is the ratio of measures between the target $(n-k)$-cell and the source $k$-cell (areas for $2$-cells, lengths for $1$-cells) for relevant geometric entities, where the measures of $0$-cells $v$ and $\tilde{v}$ are $1$ by definition. This Hodge star requires grids with circumcentric duality, which implies orthogonality. However, there has been some work to modify $\Hop$ to support non-orthogonal grids, given in \cite{Weller2014,Thuburn2014}. Unfortunately, this operator is inconsistent on many commonly used quasi-uniform spherical grids such as the cubed sphere.
    \item Barycentric Hodge star: Using a simplicial straight grid along with a barycentric twisted grid, a non-diagonal Hodge star can be defined based on nearest-neighbors. Examples of this can be found in \cite{Hirani2003,Auchmann2006}.
    \item Galerkin Hodge star \cite{Hirani2003}: Using the mass matrices induced by a set of basis functions for discrete differential forms, a Hodge star operator can be defined. This is usually done with Whitney forms \cite{Hirani2003,Tarhasaari1999}, although alternative basis functions have been developed \cite{Codecasa2016} based on energetics definitions. These approaches require a simplicial straight grid along with a barycentric twisted grid. An interesting extension of this idea to arbitrary complexes is provided in \cite{Thuburn2015}, based on primal-dual finite elements. Unfortunately, this Hodge star requires the solution of a linear system of equations, and is therefore computationally challenging for TRiSK-type schemes.
\end{itemize}

A comparison between the three types of Hodge stars for 2D surfaces can be found in \cite{Mohamed2018,Mohamed2016b}, while connections between various Hodge star definitions and how they make an inner product are explored in \cite{Hiptmair2001}.

\subsubsection{Alternative Hodge Stars} In addition to the Voronoi, barycentric and Galerkin Hodge stars, there are some additional alternative Hodge stars that have appeared in the literature. The $G_{ij}$ Hodge star \cite{Toy2017} is based on the metric tensor of the underlying coordinate system used to define the computational grid. It is defined around a topologically square, possibly non-orthogonal grid. In the case of an orthogonal grid it will reduce to the Voronoi Hodge star. The extension of this Hodge star to arbitrary topologies and geometries seems possible, although $\Hop$ will likely be quite complicated. A higher-order Hodge star on uniform rectangular grids based on coefficient fitting can be found in \cite{Deimert2016}. The extension of this idea to arbitrary topologies and geometries seems to be quite intricate, however, and it is unclear if it can be generalized. Finally, a diagonal, positive-definite Hodge star for non-orthogonal grids based on constitutive relationships can be found in \cite{ElOuafdi2021}.

\subsection{Choices for PV Wedge Products $\Rop$, $\Wop$ and $\Qop$}

The major innovation of \cite{Thuburn2009,Ringler2010}, and indeed the basis for all TRiSK-type schemes, is an explicit formula for $\Wop$ coefficients given a choice of $\Rop$ coefficients with a nearest-neighbor stencil, such that all the desirable properties for $\Wop$ are obtained ($\Wop = -\Wop^T$ and $\Rop \Dbarop{2} = \Dop{2} \Wop$). This is hereafter refered to as the TRSK2010 approach. Therefore, the choice of PV wedge product really amounts to a choice of $\Rop$, along with a choice of how to construct $\Qop$ given $\Rop$ and $\Wop$.

\subsubsection{Choice of $\Rop$/$\Wop$}
General nearest-neighbor stencil definitions of $\Rop$ and $\Wop$ are given by \cite{Hirani2003}
\begin{align}
(\Rop \twf{x}{2})_c &= \sum_{\tilde{c} \in VC(c)} R_{\tilde{c},c} \twf{x}{2}_{\tilde{c}},\\
(\Wop \twf{x}{1})_e &= \sum_{\tilde{e} \in ECP(e)} W_{\tilde{e},e} \twf{x}{1}_{\tilde{e}}. 
\end{align}
Using the TRSK2010 approach, the $\Wop$ coefficients ($W_{\tilde{e},e}$) are uniquely determined from the $\Rop$ coefficients ($R_{\tilde{c},c}$). Additionally, this choice of $\Rop$ is explicitly solvable for $\twf{q}{0}$. The specific choice of coefficients for $\Rop$ follows one of two approaches:
\begin{itemize}
    \item Metric \cite{Eldred2017,Thuburn2012}: $R_{\tilde{c},c}$ is defined in terms of overlap areas $A_{\tilde{c},c}$ between $c$ and $\tilde{c}$ and the twisted cell area $A_{\tilde{c}}$ (see Figure \ref{grid-geom-fig}) as
    \begin{eqnarray}
R_{\tilde{c},c} &=& \frac{A_{\tilde{c},c}}{A_{\tilde{c}}}.
\end{eqnarray}
    This approach is valid for arbitrary grid geometries.
    \item Combinatorical \cite{Hirani2003}: $R_{\tilde{c},c}$ is defined in terms of the size $[CV(v)]$ (number of elements) of $CV(v)$ as
    \begin{eqnarray}
    R_{\tilde{c},c} &=& \frac{1}{[CV(v)]},
    \end{eqnarray}
    where $v$ is the straight grid vertex dual to the dual grid cell $\tilde{c}$. This approach is valid for arbitrary grid geometries as well, although it has only been used, as far as we aware, on square topologies (with possible lat-lon singularities).
\end{itemize}
These two approaches give the same results on a uniform grid. Another option for $\Rop$ and $\Wop$ is to use the primal-dual FE operators from \cite{Thuburn2015}, since they share the same stencil as the DEC ones and have the same key properties.

\subsubsection{Choice of $\Qop$}
Given $\Wop$ and $\Rop$, there are four (non-exhaustive) possible definitions of $\Qop$, differing in whether they conserve total energy ($\Qop^{TE}$), potential enstrophy ($\Qop^{PE}$), both ($\Qop^{DBL}$) or none ($\Qop^{accur}$):
\begin{itemize}
    \item $\Qop^{TE}$ satifies $\Qop = -\Qop^T$, and therefore conserves total energy (as shown in Section \ref{total-energy-cons}):
  \begin{align}
(\Qop^{TE} \twf{x}{1})_e &= \sum_{\tilde{e} \in ECP(e)} \frac{q^e + q^{\tilde{e}}}{2} W_{\tilde{e},e} \twf{x}{1}_{\tilde{e}},
\end{align}
for an arbitrary dual edge reconstruction $q^e$. This definition is equivalent to
\begin{equation*}
\Qop^{TE} = \frac{1}{2} (q^e \Wop + \Wop q^e),
\end{equation*}
It is clear from this latter form of the definition that  $\Qop = -\Qop^T$, provided $\Wop = -\Wop^T$. Oftentimes, sophisticated advection schemes (APVM, LUST, CLUST, upwinding, WENO, etc.) are used to compute $q^e$ \cite{Weller2012a,Weller2012b}, especially on hexagonal grids to control the spurious Rossby mode. A clever choice of $q^e$ was developed in \cite{Gassmann2018a} to help alleviate Hollingsworth instability.
\item $\Qop^{PE}$ satisfies $\Qop \Dbarop{1} \twf{q}{0} - \Dop{1} \Rop^T \frac{(\twf{q}{0})^2}{2} = 0$,  and therefore conserves potential enstrophy (as shown in Section \ref{total-pe-cons}):
  \begin{align}
(\Qop^{PE} \twf{x}{1})_e &= \sum_{\tilde{e} \in ECP(e)} \sum_{\tilde{v} \in VE(\tilde{e})} \frac{\twf{q}{0}_{\tilde{v}}}{2} W_{\tilde{e},e} \twf{x}{1}_{\tilde{e}}.
\end{align}
\item $\Qop^{DBL}$ satisfies both $\Qop = -\Qop^T$ and $\Qop \Dbarop{1} \twf{q}{0} - \Dop{1} \Rop^T \frac{(\twf{q}{0})^2}{2}= 0$,  and therefore conserves total energy and potential enstrophy:
\begin{align}
(\Qop^{DBL} \twf{x}{1})_e &= \sum_{\tilde{e} \in ECP(e)} \sum_{\tilde{v} \in VE(\tilde{e})} Q_{\tilde{e},e,\tilde{v}} \twf{x}{1}_{\tilde{e}} \twf{q}{0}_{\tilde{v}}.
\end{align}  
 The $Q_{\tilde{e},e,\tilde{v}}$ coefficients in $\Qop^{DBL}$ are uniquely determined by the choice of $\Rop$ coefficients. The first example of this appeared in \cite{Arakawa1981}, with an extension to higher-order accuracy in the implied vorticity advection equation in \cite{Takano1984}. A general framework that combines \cite{Arakawa1981} and \cite{Takano1984} appears in \cite{Salmon2004}, which identifies other version of $\Qop$ with the same doubly conservative properties. All of these approaches are restricted to logically square topologies. An extension of \cite{Salmon2004} to arbitrary grid topologies and geometries is found in \cite{Eldred2017}. However, only metric $\Rop$ coefficients are evaluated in that work. A $\Qop$ with partial double conservation is found in \cite{Sadourny1975}: it conserves total energy but only enstrophy in the non-divergent limit, not potential enstrophy. 
\item $\Qop^{accur}$ gives up conservation in exchange for accurate advection of PV:
  \begin{align}
(\Qop^{accur} \twf{x}{1})_e &= \sum_{\tilde{e} \in ECP(e)} q^e W_{\tilde{e},e} \twf{x}{1}_{\tilde{e}}.
\end{align}
\end{itemize}
recalling $q^e$ is an arbitrary dual edge reconstruction, usually obtained from a sophisticated advection scheme such as APVM, LUST, CLUST, upwinding, WENO, etc.

For all four definitions, if $\twf{q}{0} = \twf{I}{0}$, then $\Qop$ becomes equal to the $\Wop$ defined from $\Rop$. This is clear for $\Qop^{TE}$, $\Qop^{PE}$ and $\Qop^{accur}$ since they are based on the $\Wop$ coefficients and $q^e = 1$ when $\twf{q}{0} = \twf{I}{0}$. For $\Qop^{DBL}$ see the proof in Appendix \ref{leibniz-appendix} that enforcing both total energy and potential enstrophy conservation automatically gives PV compatibility.

\subsection{Choices for KE Wedge Product}
Unlike the PV wedge product, the KE wedge product does not have to satisfy any properties. Therefore, there is some additional freedom in how it is defined. A general nearest-neighbor stencil definition of $\Wedgeop{\stf{x}{1}}{\twf{y}{1}}$ is given by
\begin{equation}
\label{KE-wedge-stencil}
(\Wedgeop{\stf{x}{1}}{\twf{y}{1}})_{\tilde{c}} = \sum_{\tilde{e} \in EC(\tilde{c})} \mathfrak{T}_{\tilde{c},e,\tilde{e}} \stf{x}{1}_e \twf{y}{1}_{\tilde{e}},
\end{equation}
where $e$ is the straight grid edge corresponding to the twisted grid edge $\tilde{e}$. As for $\Rop$, either metric or combinatorial definitions can be used for $\mathfrak{T}_{\tilde{c},e,\tilde{e}}$:
\begin{itemize}
    \item Metric \cite{Ringler2010,Thuburn2014}: $\mathfrak{T}_{\tilde{c},e,\tilde{e}}$ is defined in terms of overlap areas $A_{\tilde{c},\tilde{e}}$ and extended edge area $\tilde{A}_{\tilde{e}}$ as
    \begin{equation}
    \mathfrak{T}_{\tilde{c},e,\tilde{e}} = \frac{A_{\tilde{c},\tilde{e}}}{\tilde{A}_{\tilde{e}}}.
\end{equation}
    This approach is valid for arbitrary grid geometries.
    \item Combinatorical \cite{Hirani2003}: $\mathfrak{T}_{\tilde{c},e,\tilde{e}}$ is defined in terms of the size $[CE(\tilde{e})]$ of $CE(\tilde{e})$, which is always $\frac{1}{2}$ on a conforming grid without boundaries, i.e.
    \begin{eqnarray}
    \mathfrak{T}_{\tilde{c},e,\tilde{e}} &=& \frac{1}{2}.
    \end{eqnarray}
    This approach is valid for arbitrary grid geometries as well, although it has only been used on square topologies (with possible lat-lon singularities).
\end{itemize}
From the coefficients $\mathfrak{T}_{\tilde{c},e,\tilde{e}}$, it is easy to define the adjoint wedge products $\Wedgeop{\stf{x}{0}}{\twf{y}{1}}$ and $\Wedgeop{\stf{x}{0}}{\stf{y}{1}}$ needed to compute the mass flux $\twf{F}{n-1}$:
\begin{eqnarray}
\label{mass-flux-wedge-1}
(\Wedgeopadj{\stf{h}{0}}{\stf{u}{1}})_e &=& \sum_{\tilde{c} \in CE(e)} \mathfrak{T}_{\tilde{c},e,\tilde{e}} \stf{h}{0}_v \stf{u}{1}_e,\\
\label{mass-flux-wedge-2}
(\Wedgeopadj{\stf{h}{0}}{\twf{u}{n-1}})_{\tilde{e}} &=& \sum_{\tilde{c} \in CE(\tilde{e})} \mathfrak{T}_{\tilde{c},e,\tilde{e}} \stf{h}{0}_v \twf{u}{n-1}_{\tilde{e}}\ .
\end{eqnarray}

In the case of a uniform grid, the metric and combinatorical approaches give the same results.

\paragraph{Alternative Choices}
Clever choices of KE wedge products with slightly different stencils were developed in \cite{Hollingsworth1983,Gassmann2013} to help alleviate Hollingsworth instability. Another alternative KE wedge product with improved accuracy was developed in \cite{Engwirda2022}, based on a least square reconstruction, although it was used only to compute kinetic energy and not the corresponding mass flux. As for $\Rop$ and $\Wop$, it would also be possible to use the primal-dual FE approach from \cite{Thuburn2015} to define a KE wedge product, although this has not been worked out yet. A final view of the wedge product is through the connection with the interior product (contraction), given by Hirani's formula $\intp{\mathbf{x}} \alpha = (-1)^{k(n-k)} \tstar (\tstar \alpha \wedge \mathbf{x}^\flat)$. A discrete version of this is explored in \cite{Bossavit2003}.

\subsection{Choice of grid}
As discussed in Section \ref{grids}, grids have two aspects: topology and geometry. Within TRiSK-type schemes, there is a distinction between schemes that are capable of handling only topologically square grids (along with possibly a lat-lon grid singularity), and those capable of handling arbitrary topologies. Additionally, there is a distinction between schemes that can handle arbitrary (non-orthogonal) grid geometries, and those that can handle only orthogonal grid geometries. For quasi-uniform grids such as the icosahedral-pentagonal and cubed sphere, the ability to handle general conforming topologies (for both) and non-orthogonal geometries (for cubed sphere) is required. These restrictions on grid topology/geometry arise due to the choices made for the Hodge star and wedge product operators, as discussed in the previous subsections.

\subsection{Grid Optimization}
In addition to these choices of operators, it is also possible to apply many different optimization schemes to improve the geometry of the grid. This includes tweaking \cite{Heikes1995,Heikes2013}, centroidal voronoi tesselations \cite{Lu2012,Jacobsen2013,Du2010}, spring-dynamics \cite{Tomita2002,Iga2014}, generalized power grids \cite{engwirda2018generalised} and Hodge-optimized triangulations \cite{mullen2011hot}. Some of these optimization choices were compared in \cite{Eldred2015,Eldred2017,engwirda2018generalised}; but only for certain choices of operators. Grid optimization is primarily utilized to improve grid regularity, and therefore improve operator and scheme accuracy, and to allow a larger explicit time step.

\subsection{Effects of choices on properties}

In general, the topological properties are obtained by any of the choices discussed above, with the possible exception of total energy or potential enstrophy conservation depending on the choice of $\Qop$. The choice of operators does affect the metric properties, specifically:
\begin{itemize}
    \item The Hodge star choice (specifically $\Dop{k}$/$\Dbarop{k}$ plus $\Iop$, $\Jop$ and $\Hop$) determines the accuracy of the differential operators, and the combination of the Hodge star and the wedge product (specifically  $\Rop$/$\Wop$/$\Qop$ or $\Wedgeop{\stf{x}{1}}{\twf{y}{1}}$ plus $\Iop$, $\Jop$ and $\Hop$) determines the accuracy of the scalar/vector product operators.
    \item Hollingsworth instability is controlled by the choice of wedge products used in the kinetic energy/mass flux ($\Wedgeop{\stf{x}{1}}{\twf{y}{1}}$) and/or the PV flux term ($\Qop$), as well as Hodge stars $\Hop$ and $\Iop$.
    \item The choice of grid topology determines the presence or absence of spurious branches of dispersion relationship.
    \item The choice of grid optimization can affect the accuracy of the differential operators.
\end{itemize}

The main issues seen with TRiSK-type schemes are the inability to find (on arbitrary grids) accurate and Hollingsworth free discrete versions of the operators. However, as discussed in Section \ref{unexplored}, only a small portion of the possible choices have been explored. With careful grid optimization it seems possible to, at least, get 1st order in $L_\infty$ for some tested quasi-uniform spherical grids for differential operators, but scalar/vector products are trickier. The combinatorial wedge products do not have accuracy issues, including lat-lon variants and the extension \cite{Toy2017} to topologically square non-orthogonal grids. However, these wedge products have not been tested on arbitrary grids. The metric wedge products ($\Rop$/$\Wop$/$\Qop$ and $\Wedgeop{\stf{x}{1}}{\twf{y}{1}}$) on arbitrary grids do have accuracy issues in at least $L_\infty$. Hollingsworth instability appears avoidable (or at least controllable) through careful specification of either $q^e$ using in $\Qop^{TE}$ \cite{Gassmann2018a} and/or $\Wedgeop{\stf{x}{1}}{\twf{y}{1}}$ \cite{Hollingsworth1983,Gassmann2013}.

\section{Schemes from the literature}
\label{literature-schemes}
In this section, we identify the operator and grid optimization choices that were made\footnote{Usually implicitly, unaware of the DEC structure.} in a variety of TRiSK-type schemes from the literature. Specifically, we will look at the schemes in 
\cite{Arakawa1981} (AL81), \cite{Sadourny1975} (SD75), \cite{Takano1984} (TW84), \cite{Ringler2010} (TRSK2010), \cite{Thuburn2014} (Thuburn2014) , \cite{Weller2012a,Weller2012b,Weller2014} (Weller2014), \cite{Eldred2015,Eldred2017} (Eldred2017) \cite{Salmon2004} (S04) and \cite{Toy2017} (Toy2017). Some of this connection was covered in \cite{Eldred2015,Eldred2017}, but the understanding of the discrete wedge product was lacking. A detailed summary of the various choices made by TRiSK-type schemes in the literature is given in Table \ref{scheme-summary-table}. 

\begin{table}
\begin{center}
\begin{tabular}{ |c|c|c|c|c|c|c|c| } 
 \hline
Scheme & Grid & $\Iop$/$\Jop$/$\Hop$ & $\Rop$/$\Wop$ & $\Qop$ & $\Wedgeop{\stf{x}{1}}{\twf{y}{1}}$ & Dof Scaling & Grid Opt\\ \hline \hline
 AL81 & S-O/LL & V & C & DBL & C & yes & N/A\\ \hline
 SD75 & S-O/LL & V & C & DBL$^1$ & C & yes & N/A\\ \hline
 TW84 & S-O/LL & V & C & DBL$^3$ & C & yes & N/A\\ \hline
 TRSK2010 & US-O & V & M & TE/PE & M & yes & CVT \\ \hline
 Thuburn2014 & US-NO & V$^2$ & M & ACCUR$^4$ & other$^4$ & no & CVT\\ \hline
 Weller2014 & US-NO & V$^2$ & M & TE & other$^5$ & partial$^5$ & CVT\\ \hline
 Eldred2017 & US-NO & V$^2$ & M & DBL/TE/PE & M & no & CVT, SD, TW\\ \hline
 S04 & S-O & V & C & DBL$^3$ & C & yes & N/A\\ \hline
 Toy2017 & S-NO & $G_{ij}$ & C & DBL & C & yes & N/A\\
 \hline
\end{tabular}
\end{center}
\caption{A summary of the TRiSK-type schemes in the literature: the grid topologies and geometries they support; the choice of Hodge stars $\Iop$/$\Jop$/$\Hop$, PV wedge products $\Rop$/$\Wop$ and $\Qop$, KE wedge product $\Wedgeop{\stf{x}{1}}{\twf{y}{1}}$ they make; whether dof scaling is required to put them into the DEC framework used in this paper; and the choice of grid optimization. Here S-O indicates square topology with orthogonal geometry, LL indicates lat-lon topology and geometry, US-O indicates an unstructured topology with orthogonal geometry, US-NO indicates an unstructured topology with non-orthogonal geometry and S-NO indicates a square topology with non-orthogonal geometry. A V indicates a Voronoi Hodge star, while $G_{ij}$ indicates the $G_{ij}$ based Hodge star. A C indicates a combinatorial discretization for the PV or KE wedge product, while an M indicates a metric one. For $\Qop$, DBL indicates the doubly conservative version $\Qop^{DBL}$, TE indicates the total energy version $\Qop^{TE}$, PE indicates the potential enstrophy conserving version $\Qop^{PE}$ and ACCUR indicates the accurate PV advection version $\Qop^{accur}$. CVT indicates centroidal Voronoi tessellation grid optimization, SD indicates spring dynamics grid optimization and TW indicates the tweaking grid optimization. \\
$^1$ The SD75 scheme only conserves total energy plus enstrophy in the divergent limit, not potential enstrophy.\\
$^2$ The Thuburn2014 and Eldred2017 schemes use a modified $\Hop$ for non-orthogonal grids, that reduces to the diagonal Voronoi $\Hop$ for an orthogonal grid. Weller2014 explores the same modified $\Hop$, but also introduces non-symmetric $\Hop$ choices that do not conserve energy.\\
$^3$ The TW84 and S04 schemes are 4th-order for the implied vorticity equation (only optionally in the case of S04).\\
$^4$ The Thuburn2014 scheme gives up energy conservation in exchange for accurate nonlinear advection (through choice of $\Qop$ and mass flux), and computes kinetic energy using a least-squares approach.\\
$^5$ The Weller2014 scheme computes kinetic energy using a metric wedge product, but uses an accurate advection scheme for mass flux $\twf{F}{n-1}$; and requires scaling of fluid height but not velocity to fit into the DEC framework.
}
\label{scheme-summary-table}
\end{table}

\subsection{Degree of freedom scaling}
Many TRiSK-type schemes in the literature work with point values, instead of the integral quantities $\twf{h}{2}$ and $\stf{v}{1}$. Therefore, for these schemes to be put into the same form as the general DEC framework the degrees of freedom must be scaled. This involves multiplying the predicted pointwise height by the area of a twisted $2$-cell $A_{\tilde{c}}$ to obtain $\twf{h}{2}$, and the predicted pointwise velocity by the length of a straight $1$-cell (edge) $A_e$ to obtain $\stf{v}{1}$.

\subsection{Specific schemes}

\paragraph{AL81} The AL81 scheme \cite{Arakawa1981} uses a Voronoi Hodge star, combinatorical PV and KE wedge products and $\Qop^{DBL}$. It was developed for square, orthogonal grids, including the modifications needed to handle polar singularities on a latitude-longitude grid. A modification of the KE wedge product to suppress Hollingsworth instability can be found in \cite{Hollingsworth1983,Lazic1986}.

\paragraph{SD75} The SD75 scheme \cite{Sadourny1975} uses a Voronoi Hodge star, combinatorical PV and KE wedge products and a slightly modified $\Qop^{DBL}$ that only conserves total energy and enstrophy in the non-divergent limit, not potential enstrophy. It was developed for square, orthogonal grids, including the modifications needed to handle polar singularities on a latitute-longitude grid.

\paragraph{TW84} The TW84 scheme \cite{Takano1984} extended the AL81 scheme with a  $\Qop^{DBL}$ that is 4th order for the implied vorticity equation. It is otherwise the same, including the choice of grids.

\paragraph{S04} The S04 scheme \cite{Salmon2004} generalizes AL81 and TW84 to a family of $\Qop^{DBL}$ (optionally extended to 4th order for the implied vorticity equation) operators, through a quasi-Hamiltonian approach. It therefore also uses a Voronoi Hodge star and combinatorical PV and KE wedge products. It was developed for square, orthogonal grids.

\paragraph{TRSK2010} The TRSK2010 scheme \cite{Ringler2010} uses a Voronoi Hodge star, metric PV and KE wedge products, and either $\Qop^{TE}$ or $\Qop^{PE}$. It is designed for unstructured, orthogonal grids, primarily Delauney triangulations and their Voronoi duals such as the icosahedral-hexagonal grid. These grids are usually optimized using the centroidal Voronoi tesselation approach. TRSK2010 can be viewed as an extension of AL81 to unstructured grids with only partial conservation. For $\Qop^{TE}$, a variety of choices for $q^e$ were explored, focusing mainly on the anticipated potential vorticity method (APVM). A modification of the KE wedge product to suppress Hollingsworth instability \cite{Gassmann2013}, and clever choice of $q^e$ to do the same for $\Qop^{TE}$ can be found in \cite{Gassmann2018a}.
 
\paragraph{Thuburn2014} The Thuburn2014 scheme \cite{Thuburn2014} uses a Voronoi Hodge star (with a possible modification of $\Hop$ to support non-orthogonal grids) and metric PV wedge products. However, it gives up conservation of total energy and potential enstrophy for a more accurate computation of kinetic energy and better nonlinear advection properties. Specifically, the kinetic energy is computed using a least-squares approach; and accurate advection schemes are used to compute $q^e$ and the mass flux, using $\Qop^{accur}$. It is designed for unstructured, non-orthogonal grids, including the icosahedral-hexagonal and cubed-sphere grids. The grid is optimized using the centroidal Voronoi tesselation approach or a related approach on cubed-sphere grids. This scheme can be thought of as an extension of TRSK2010 to non-orthogonal grids that gives up energy conservation for accurate kinetic energy and better advection.

\paragraph{Weller2014} The Weller2014 scheme \cite{Weller2012a,Weller2012b,Weller2014} uses metric PV wedge products and $\Qop^{TE}$. Voronoi Hodge stars are used for $\Iop$ and $\Jop$, and several options are explored for $\Hop$: the modified $\Hop$ from \cite{Thuburn2014}, and a non-symmetric $\Hop$ that reduces to the Voronoi $\Hop$ in the case of an orthogonal grid. A non-symmetric $\Hop$ leads to a loss of energy conservation. The kinetic energy is computed using a metric wedge product, while the mass flux $\twf{F}{n-1}$ is computed using an accurate advection scheme. This mass flux also gives up energy conservation. In $\Qop^{TE}$, $q^e$ is computed for a variety of different advection schemes. It is designed for unstructured, non-orthogonal grids, including the icosahedral-hexagonal and cubed-sphere grids. The grid is optimized used the centroidal Voronoi tesselation approach or a related approach on cubed-sphere grids. This scheme can be thought of as an extension of TRSK2010 to non-orthogonal grids that gives up energy conservation for better advection.

\paragraph{Eldred2017} The Eldred2017 scheme \cite{Eldred2015,Eldred2017} uses a Voronoi Hodge star (with a possible modification of $\Hop$ to support non-orthogonal grids from \cite{Thuburn2014}), metric PV and KE wedge products and any of $\Qop^{TE}$, $\Qop^{PE}$ or $\Qop^{DBL}$. It is essentially a merger of S04 and TRSK2010, with the main novelty being an extension of $\Qop^{DBL}$ to arbitrary grids. It is designed for unstructured, non-orthogonal grids, including the icosahedral-hexagonal and cubed-sphere grids. The grid is optimized using the centroidal Voronoi tesselation, spring dynamics or tweaking approaches.

\paragraph{Toy2017} The Toy2017 scheme \cite{Toy2017} uses a $G_{ij}$ Hodge star, combinatorical PV and KE wedge products and $\Qop^{DBL}$. It is designed for a square, non-orthogonal grid; and primarily tested for smooth deformations of a uniform square grid. The new Hodge star reduces to a Voronoi Hodge star for orthogonal grids, and therefore the scheme can be thought of as an extension of AL81 to non-orthogonal grids. 

\subsection{Unexplored approaches}
\label{unexplored}
This (brief) review of existing schemes immediately highlights some unexplored operator choices that fit within the DEC framework. Specifically, for the Hodge stars:
\begin{itemize}
    \item barycentric Hodge stars \cite{Codecasa2016} (straightforward, already developed);
    \item $G_{ij}$ Hodge star from \cite{Toy2017} on arbitrary grids (tricky, not yet developed);
    \item constitutive relationship based Hodge star from \cite{ElOuafdi2021} (straightforward, already developed);
\end{itemize}
for the PV wedge product (specifically $\Rop$ and $\Wop$):
\begin{itemize}
    \item combinatorical wedge products on arbitrary grids (straightforward, already developed),
    \item primal-dual FE wedge product operators from \cite{Thuburn2015} (straightforward, already developed);
    \end{itemize}
and for the KE wedge product:
\begin{itemize}
    \item combinatorical wedge products on arbitrary grids (straightforward, already developed),
    \item primal-dual FE wedge product operators from \cite{Thuburn2015} (straightforward, not yet developed),
    \item least squares reconstruction wedge product from \cite{Engwirda2022} (straightforward, already developed),
    \item extrusion/contraction based wedge product from \cite{Bossavit2003} (straightforward, already developed).
\end{itemize}
Most of these options have already been developed, just not studied in the context of TRiSK-type schemes. Only the primal-dual FE versions of the KE wedge product (straightforward) and the $G_{ij}$ Hodge star on arbitrary grids (tricky) would require further research. These alternative choices do not expand the stencil of the general scheme beyond nearest-neighbor (although the Hodge star is no longer diagonal) and keep the topological properties of the scheme. 

These alternative choices of operators are interesting for several reasons:
\begin{itemize}
    \item Accuracy results suggest that the combinatorical approach to wedge products might be superior, although this is yet to be verified for quasi-uniform spherical grids. Specifically, \cite{Toy2017} uses combinatorical wedge products and is accurate on highly-distorted non-orthogonal planar grids, while \cite{Eldred2017} showed that metric wedge products are inaccurate on quasi-uniform spherical grids.
    A combinatorical definition also fits with differential geometry: the wedge product is topological.
    \item Barycentric Hodge stars expand the possible grid geometries, which is especially important for grids with complicated lateral boundaries or highly deformed grid cells.
    \item The primal-dual FE scheme operators have proven accuracy \cite{Thuburn2015}, albeit when used in scheme with mass matrices. 
    \item The least squares based wedge product has proven accuracy \cite{Engwirda2022}, albeit when used in a non-mimetic scheme.
\end{itemize}
Therefore, it seems possible that some combination of these choices for operators and grid optimizations could obtain acceptable accuracy on arbitrary grids and hopefully avoid the Hollingsworth instability.

\section{Conclusions}
\label{conclusions}
 In this paper, we have shown that TRiSK-type schemes are best understood as a discrete exterior calculus applied to a Hamiltonian formulation based on split exterior calculus. The key to this was the wedge product and the topological pairing, with the latter newly introduced here in this manuscript. This new understanding gives a clear separation between the topological Poisson brackets, written entirely in terms of the wedge product, exterior derivative and topological pairing; and the metric Hamiltonian, which contains, in addition, the Hodge star and inner product. This separation is somewhat mimicked in the discrete setting, through the use of combinatorial operators for the topological parts and metric operators for the metric parts. However, the separation can be broken in the case of the wedge product, where metric operators are sometimes used. All known TRiSK-type schemes can be written in terms of this framework, differing only through their choice of Hodge stars and wedge products, along with the grids supported by these choices. 
 
 This new framework has several immediate possible applications. The first is a comprehensive investigation of accuracy and Hollingsworth instability for both explored and unexplored combinations from Section \ref{unexplored}, on a variety of grids with various optimization choices. The main goal would be to identify a combination of grid choice/optimization, wedge products and Hodge stars that gives acceptable accuracy and avoids Hollingsworth instability. Such a scheme could serve as a revised basis for the operational atmospheric and oceanic models using TRiSK-type schemes that offers improved accuracy while keeping the desirable properties and nearest-neighbor stencil (and thus computational efficiency).
 
A second application is exploring TRiSK-type schemes with lateral boundaries, which are needed for ocean models and limited-area atmospheric models. Existing treatment of boundaries in DEC (and therefore TRiSK-type schemes) is somewhat ad-hoc, and lacks the generality to treat arbitrary boundary conditions in a consistent way. In recent work \cite{EldredSandBnd}, a DEC with consistent treatment of boundaries and arbitrary boundary conditions was developed. The application of this DEC to TRiSK-type schemes will be reported in a future publication.

A final application is studying whether the (nonlinear) conservation properties of total energy and potential enstrophy can be kept while introducing accurate advection schemes for the fluid height and the potential vorticity, or in other words, can the loss of conservation properties in Thuburn2014 and Weller2014 be remedied in a way that still keeps accurate advection? A positive answer to this is given in \cite{EldredCRMReport}, where a variant of TRiSK with WENO/FCT advection (that also has high-order Hodge stars) is developed for uniform rectangular grids. This does require expanding the stencil of the advection operators and the Hodge stars, but no more than required by the advection operators in the Weller2014 or Thuburn2014 schemes. Work is currently ongoing to extend this idea to unstructured grids.

\section{Acknowledgements}
This article has been co-authored by an employee of National Technology \& Engineering Solutions of Sandia, LLC under Contract No. DE-NA0003525 with the U.S. Department of Energy (DOE). The employee owns right, title and interest in and to the article and is responsible for its contents. The United States Government retains and the publisher, by accepting the article for publication, acknowledges that the United States Government retains a non-exclusive, paid-up, irrevocable, world-wide license to publish or reproduce the published form of this article or allow others to do so, for United States Government purposes. The DOE will provide public access to these results of federally sponsored research in accordance with the DOE Public Access Plan https://www.energy.gov/downloads/doe-public-access-plan.

\bibliographystyle{abbrv}
\bibliography{main,main1}

\appendix
    
\section{Proof of $\Dop{k}$/$\Dbarop{k}$/$\topopair{}{}$ Properties}
\label{diff-properties}
Start by noting that the co-boundary operator versions of the discrete exterior derivatives $\Dop{k}$ and $\Dbarop{k}$ satisfy by construction \cite{Hirani2003}:
\begin{equation}
\label{D-transpose}
    \Dbarop{2} = -\Dop{1}^T , \quad\quad\quad \Dop{2} = \Dbarop{1}^T ,
\end{equation}
and
    \begin{equation}
    \label{D-const2}
        \Dop{1} \stf{c}{0} = 0 , \quad\quad\quad \Dbarop{1} \twf{c}{0} = 0 ,
    \end{equation}
for arbitrary constants $\stf{c}{0}$ and $\twf{c}{0}$, which are just (\ref{D-const}) and (\ref{D-adjoints}). Using (\ref{D-transpose}) plus the definition of the topological pairing (\ref{discrete-topo-pairing-1}) - (\ref{discrete-topo-pairing-2}), it is easy to show that integration by parts (\ref{D-IBP1}) and (\ref{D-IBP2}) hold:
    \begin{eqnarray}
    \topopair{\Dop{k} \stf{x}{k-1}}{\twf{y}{n-k}} + (-1)^{k-1} \topopair{\stf{x}{k-1}}{\Dbarop{n-k+1} \twf{y}{n-k}} = 0, \\
    \topopair{\Dbarop{k} \twf{x}{k-1}}{\stf{y}{n-k}} + (-1)^{k-1} \topopair{\twf{x}{k-1}}{\Dop{n-k+1} \stf{y}{n-k}} = 0.
    \end{eqnarray}
Finally, applying integration by parts for the case $k=2$ (recalling we are in $n=2$) plus (\ref{D-const2}) with $\twf{y}{n-k} = \twf{I}{0}$ and $\stf{y}{n-k} = \stf{I}{0}$ gives the discrete Stokes theorems (\ref{D-stokes1}):
    \begin{equation}
(\twf{I}{0})^T \Dop{2} \stf{x}{1} = 0 , \quad\quad
(\stf{I}{0})^T \Dbarop{2} \twf{x}{1} = 0 .
\end{equation}

\section{Relationship between PV wedge product Leibniz rules}
\label{leibniz-appendix}
Consider a discrete wedge product $\Wedgeop{\twf{x}{0}}{\twf{y}{1}}$ that satisfies the full Leibniz rule
\begin{equation}
\label{full-leibniz}
    \Wedgeop{\twf{x}{0}}{\Dbarop{1} \twf{y}{0}} + \Wedgeop{\twf{y}{0}}{\Dbarop{1} \twf{x}{0}} = \Dop{1} (\Wedgeopadj{\twf{x}{0}}{\twf{y}{0}}) 
\end{equation}
along with anti-symmetry (\ref{Q-antisymmetry}), where $\Wedgeopadj{\twf{x}{0}}{\twf{y}{0}}$ is the adjoint of $\Wedgeop{\twf{x}{0}}{\twf{y}{2}}$. Given coefficients for $\Wedgeopadj{\twf{x}{0}}{\twf{y}{0}}$ (i.e. really coefficients for $\Wedgeop{\twf{x}{0}}{\twf{y}{2}}$), it seems likely that $\Wedgeop{\twf{x}{0}}{\twf{y}{1}}$ can be defined in such a way that (\ref{full-leibniz}) and (\ref{Q-antisymmetry}) are both satisfied, following the same sort of general procedure as in \cite{Eldred2017}. In fact, it is hypothesized that the $\Wedgeop{\twf{x}{0}}{\twf{y}{1}}$ wedge product defined in \cite{Eldred2017} satisfies the more general conditions (\ref{full-leibniz}) rather than just (\ref{Q-pens}). Since, as shown below, a wedge product that satisfies (\ref{full-leibniz}) also satisfies both (\ref{Q-pvcompat}) and (\ref{Q-pens}), this would explain why PV compatibility (\ref{Q-pvcompat}) did not have to be enforced in \cite{Eldred2017}, but was found to hold anyways.

\subsection{Partial Leibniz rule (\ref{Q-pvcompat})}
Start by taking the topological pairing of (\ref{full-leibniz}) with an arbitrary twisted $1$-form $\twf{z}{1}$:
\begin{equation}
\topopair{\twf{z}{1}}{\Wedgeop{\twf{x}{0}}{\Dbarop{1} \twf{y}{0}}} + \topopair{\twf{z}{1}}{\Wedgeop{\twf{y}{0}}{\Dbarop{1} \twf{x}{0}}} = \topopair{\twf{z}{1}}{\Dop{1} (\Wedgeopadj{\twf{x}{0}}{\twf{y}{0}})}.
\end{equation}
Now let $\twf{x}{0} = \twf{I}{0}$ such that
\begin{equation}
\topopair{\twf{z}{1}}{\Wedgeop{\twf{I}{0}}{\Dbarop{1} \twf{y}{0}}} + \topopair{\twf{z}{1}}{\Wedgeop{\twf{y}{0}}{\Dbarop{1} \twf{I}{0}}} = \topopair{\twf{z}{1}}{\Dop{1} (\Wedgeopadj{\twf{I}{0}}{\twf{y}{0}})},
\end{equation}
and use (\ref{D-const2}) to eliminate the middle term yielding
\begin{equation}
\topopair{\twf{z}{1}}{\Wedgeop{\twf{I}{0}}{\Dbarop{1} \twf{y}{0}}} = \topopair{\twf{z}{1}}{\Dop{1} (\Wedgeopadj{\twf{I}{0}}{\twf{y}{0}})}.
\end{equation}
Apply IBP (\ref{D-IBP1}) and (\ref{D-IBP2}) on the rhs which gives
\begin{equation}
\topopair{\twf{z}{1}}{\Wedgeop{\twf{I}{0}}{\Dbarop{1} \twf{y}{0}}} = \topopair{\Dbarop{2} \twf{z}{1}}{ \Wedgeopadj{\twf{I}{0}}{\twf{y}{0}}}.
\end{equation}
Now use the definition of operator adjoints (\ref{discrete-adjoint}) or \eqref{adjoint-2} on the right term and antisymmetry (\ref{Q-antisymmetry}) on the lhs to get
\begin{equation}
-\topopair{\Dbarop{1} \twf{y}{0} }{\Wedgeopadj{\twf{I}{0}}{\twf{z}{1}}} = \topopair{\twf{y}{0}}{ \Wedgeop{\twf{I}{0}}{\Dbarop{2} \twf{z}{1}}}.
\end{equation}
Finally, apply IBP (\ref{D-IBP1}) and (\ref{D-IBP2}) on the lhs to result in
\begin{equation}
\topopair{\twf{y}{0}}{\Dop{2} (\Wedgeopadj{\twf{I}{0}}{\twf{z}{1}})} = \topopair{\twf{y}{0}}{\Wedgeop{\twf{I}{0}}{\Dbarop{2} \twf{z}{1}}}.
\end{equation}
Since this must hold for arbitrary $\twf{y}{0}$, we get
\begin{equation}
\Dop{2} (\Wedgeopadj{\twf{I}{0}}{\twf{z}{1}}) = \Wedgeop{\twf{I}{0}}{\Dbarop{2} \twf{z}{1}},
\end{equation}
which is just (\ref{Q-pvcompat}).

\subsection{Partial Leibniz rule (\ref{Q-pens})}
The partial Leibniz rule (\ref{Q-pens}) is just a special case of (\ref{full-leibniz}) with $\twf{x}{0} = \twf{y}{0} = \twf{q}{0}$:
\begin{equation}
    \Wedgeop{\twf{q}{0}}{\Dbarop{1} \twf{q}{0}} + \Wedgeop{\twf{q}{0}}{\Dbarop{1} \twf{q}{0}} = \Dop{1} (\Wedgeopadj{\twf{q}{0}}{\twf{q}{0}}) ,
\end{equation}
which by inspection is equal to (\ref{Q-pens}).

\section{Derivation of $\Hh$ functional derivatives}
\label{functional-derivs-appendix}
Start by splitting the Hamiltonian into two parts: a kinetic energy and a potential energy as
\begin{equation}
    \Hh[\stf{v}{1}, \twf{h}{2}] = \Hh_p[\twf{h}{2}]  + \Hh_{kin}[\stf{v}{1}, \twf{h}{2}] ,
\end{equation}
where
\begin{equation}
    \Hh_{kin}[\stf{v}{1}, \twf{h}{2}] = \innerprod{\twf{h}{2}}{\frac{\Wedgeop{\stf{u}{1}}{\twf{u}{n-1}}}{2}} , \quad\quad\quad \Hh_p[\twf{h}{2}] = \frac{g}{2} \innerprod{\twf{h}{2}}{\twf{h}{2}} +  g \innerprod{\twf{h}{2}}{\twf{h_s}{2}},
\end{equation}
recalling that $\stf{v}{1} = \stf{u}{1}+ \stf{R}{1}$. The variations of $\Hh_p$ are simple
\begin{equation}
    \delta \Hh_p = g \innerprod{\delta \twf{h}{2}}{\twf{h}{2}} +  g \innerprod{\delta \twf{h}{2}}{\twf{h_s}{2}},
\end{equation}
and yield
\begin{equation}
    \stf{B}{0}_p := \twdede{\Hh_p}{\twf{h}{2}} = g (\stf{h}{0} + \stf{h_s}{0}) .
\end{equation}
The variations of $\Hh_{kin}$ are more complicated:
\begin{equation}
\delta \Hh_{kin} = \innerprod{\delta \twf{h}{2}}{\frac{\Wedgeop{\stf{u}{1}}{\twf{u}{n-1}}}{2}} + \innerprod{\twf{h}{2}}{\frac{\Wedgeop{\delta \stf{u}{1}}{\twf{u}{n-1}}}{2}} + \innerprod{\twf{h}{2}}{\frac{\Wedgeop{\stf{u}{1}}{\delta \twf{u}{n-1}}}{2}}.
\end{equation}
The first term gives
\begin{equation}
    \stf{B}{0}_{kin} := \twdede{\Hh_{kin}}{\twf{h}{2}} = \frac{1}{2} \Hodgebarop{2} (\Wedgeop{\stf{u}{1}}{\twf{u}{n-1}});
\end{equation}
the remaining terms are
\begin{equation}
    \innerprod{\twf{h}{2}}{\frac{\Wedgeop{\delta \stf{u}{1}}{\twf{u}{n-1}}}{2}} + \innerprod{\twf{h}{2}}{\frac{\Wedgeop{\stf{u}{1}}{\delta \twf{u}{n-1}}}{2}}.
\end{equation}
We continue by writing these in terms of topological pairings as
\begin{equation}
    \topopair{\stf{h}{0}}{\frac{\Wedgeop{\delta \stf{u}{1}}{\twf{u}{n-1}}}{2}} + \topopair{\stf{h}{0}}{\frac{\Wedgeop{\stf{u}{1}}{\delta \twf{u}{n-1}}}{2}}.
\end{equation}
Now use the definition of operator adjoints (\ref{discrete-adjoint}) and \eqref{adjoint-2} to get
\begin{equation}
    \topopair{\delta \stf{u}{1}}{\frac{\Wedgeopadj{\stf{h}{0}}{\twf{u}{n-1}}}{2}} + \topopair{\delta \twf{u}{n-1}}{\frac{\Wedgeopadj{\stf{u}{1}}{\stf{h}{0}}}{2}}.
\end{equation}
Since $\delta \twf{u}{n-1} = \Hodgeop{1} \delta \stf{u}{1}$ and $\twf{u}{n-1} = \Hodgeop{1} \stf{u}{1}$, this gives
\begin{equation}
   \twf{F}{n-1}_{kin} = \twdede{\Hh_{kin}}{\stf{v}{1}} = \frac{1}{2} \Wedgeopadj{\stf{h}{0}}{\Hodgeop{1} \stf{u}{1}} + \frac{1}{2} \Hodgeop{1} (\Wedgeopadj{\stf{h}{0}}{\stf{u}{1}}).
\end{equation}
Finally, we can obtain $\twf{F}{n-1}$ and $\stf{B}{0}$ as
\begin{equation}
     \twf{F}{n-1} := \twdede{\Hh}{\stf{v}{1}} = \twf{F}{n-1}_{kin} ,
\end{equation}
and
\begin{equation}
    \stf{B}{0} := \twdede{\Hh}{\twf{h}{2}} = \stf{B}{0}_p + \stf{B}{0}_{kin} ,
\end{equation}
which are the same as (\ref{discrete-func-derivs}).

For alternative choice of $\Hh_{kin}$ similar procedures would yield the same form for the functional derivatives, just with a change of which wedge product is primary ($\Wedgeop{}{}$) and which ones are adjoints ($\Wedgeopadj{}{}$).

\subsection{Simplified form for $\twdede{\Hh}{\stf{v}{1}}$}
\label{simplified-form}
In the continuous RSWE, the mass flux $\frac{1}{2} \left[ \stf{h}{0} \wedge \twf{u}{n-1} + \tstar (\stf{h}{0} \wedge \stf{u}{1}) \right]$ can be simplified as $\stf{h}{0} \wedge \twf{u}{n-1}$. A discrete version of this would require that
\begin{equation}
\label{discrete-mass-flux-collapse}
    \Hodgeop{1} (\Wedgeopadj{\stf{h}{0}}{\stf{u}{1}})  = \Wedgeopadj{\stf{h}{0}}{\Hodgeop{1} \stf{u}{1}},
\end{equation}
and therefore $\twf{F}{n-1} = \Wedgeopadj{\stf{h}{0}}{\Hodgeop{1} \stf{u}{1}} = \Wedgeopadj{\stf{h}{0}}{\twf{u}{n-1}}$. To determine the conditions on the operators required, recall the adjoint wedge product for the mass flux (\ref{mass-flux-wedge-1}) - (\ref{mass-flux-wedge-2}) defined from the nearest-neighbor stencil form of KE wedge product (\ref{KE-wedge-stencil}):
\begin{eqnarray}
\label{mass-flux-adjoint-app-1}
(\Wedgeopadj{\stf{h}{0}}{\stf{u}{1}})_e &=& \sum_{\tilde{c} \in CE(e)} \mathfrak{T}_{\tilde{c},e,\tilde{e}}  \stf{h}{0}_v \stf{u}{1}_e,\\
\label{mass-flux-adjoint-app-2}
(\Wedgeopadj{\stf{h}{0}}{\twf{u}{n-1}})_{\tilde{e}} &=& \sum_{\tilde{c} \in CE(\tilde{e})} \mathfrak{T}_{\tilde{c},e,\tilde{e}}  \stf{h}{0}_v \twf{u}{n-1}_{\tilde{e}} ,
\end{eqnarray}
where $v$ is the straight grid vertex corresponding to twisted grid cell $\tilde{c}$. Inserting (\ref{mass-flux-adjoint-app-1}) - (\ref{mass-flux-adjoint-app-2}) plus the definition of $\Hodgeop{1}$ into (\ref{discrete-mass-flux-collapse}) gives an equation \begin{equation}
    \sum_{e} H_{\tilde{e}, e} \sum_{\tilde{c} \in CE(e)} \stf{h}{0}_v \stf{u}{1}_e = \sum_{\tilde{c} \in CE(\tilde{e})} \stf{h}{0}_v \sum_{e^\prime} H_{\tilde{e}, e^\prime} \stf{u}{1}_{e^\prime}
\end{equation}
that must hold at each twisted grid edge $\tilde{e}$. The two sides of this equation will be equal only for a diagonal Hodge star (such as the Voronoi Hodge star), otherwise there will be unmatched $\stf{h}{0}_v$ terms on the left hand side. The conditions for a more general definition of KE wedge product could be derived in a similar way, but it seems quite unlikely that anything but a diagonal Hodge star would ensure (\ref{discrete-mass-flux-collapse}) holds.

\section{2D Split Exterior $\leftrightarrow$ Vector Calculus Identities}
\label{split-vector-identities}
The relationships between 2D split exterior calculus and 2D vector calculus are collected here for reference.  Start by identifying a scalar $f$ with a straight $0$-form $\stf{f}{0}$ and a pseudoscalar $\tilde{f}$ with a twisted $0$-form $\twf{f}{0}$. As discussed in Section \ref{vector-proxies}, given the twisted volume form $\vform$, the vector $\mathbf{x}$, and the pseudovector $\tilde{\mathbf{x}}$, there are four vector proxies: the straight $1$-form $\stf{x}{1} = \mathbf{x}^\flat$ and twisted $n-1$ form $\twf{x}{n-1} = \intp{\mathbf{x}} \vform = \tstar \stf{x}{1}$ associated with $\mathbf{x}$; and the twisted $1$-form $\twf{x}{1} = \mathbf{\tilde{x}}^\flat$ and straight $n-1$ form $\stf{x}{n-1} = \intp{\mathbf{\tilde{x}}} \vform = - \tstar \twf{x}{1}$ associated with $\tilde{\mathbf{x}}$. In 2D, $n-1=1$ and $1$-form and $(n-1)$-forms appear to be the same object, leading to much confusion. Therefore we retain the notation $\stf{x}{n-1}$ and $\twf{x}{n-1}$ to properly distinguish between $1$ and $(n-1)$-forms. We can convert between twisted and straight forms using $\twf{I}{0} \wedge \stf{x}{k} = \twf{x}{k}$ and $\twf{I}{0} \wedge \twf{x}{k} = \stf{x}{k}$, and between straight 0-forms and twisted $n$-forms using $\twf{x}{n} = \stf{x}{0} \wedge \vform = \tstar \stf{x}{0}$.

In 2D (i.e. $n=2$) the following relationships hold (see \cite{Abraham2012} for proofs), where $f$ and $g$ are scalars, $\tilde{f}$ and $\tilde{g}$ are pseudoscalars, $\mathbf{x}$ and $\mathbf{y}$ are vectors and $\tilde{\mathbf{x}}$ and $\mathbf{\tilde{y}}$ are pseudovectors:
\begin{itemize}
    \item gradient, skew-gradient, divergence and curl in vector calculus $\leftrightarrow$ Hodge star and exterior derivative in split exterior calculus:
\begin{align}
& (\nabla f)^\flat = \diff \stf{f}{0}, \quad\quad
(\nabla \tilde{f})^\flat = \diff \twf{f}{0}, \quad\quad
(\nabla^\perp f)^\flat = \tstar \diff \stf{f}{0}, \quad\quad
(\nabla^\perp \tilde{f})^\flat = \tstar \diff \twf{f}{0} , \\
& \nabla^\perp \cdot \mathbf{x} = \tstar \diff \stf{x}{1} , \quad\quad
\nabla^\perp \cdot \mathbf{\tilde{x}} = \tstar \diff \twf{x}{1} , \\
& \nabla \cdot \mathbf{x} = \tstar \diff \tstar \stf{x}{1} = \tstar \diff \twf{x}{n-1} , \quad\quad
\nabla \cdot \mathbf{\tilde{x}} = \tstar \diff \tstar \twf{x}{1} = \tstar \diff \stf{x}{n-1} ;
\end{align}
\item scalar product in vector calculus $\leftrightarrow$ Hodge star and wedge product in split exterior calculus:
\begin{equation}
f g = \stf{f}{0} \wedge \stf{g}{0} , \quad\quad
\tilde{f} g = \twf{f}{0} \wedge \stf{g}{0}, \quad\quad
\tilde{f} \tilde{g} = \twf{f}{0} \wedge \twf{g}{0}; 
\end{equation}
\item scalar-vector product in vector calculus $\leftrightarrow$ Hodge star and wedge product in split exterior calculus:
\begin{align}
&f \mathbf{x} = \stf{f}{0} \wedge \stf{x}{1}, \quad\quad
\tilde{f} \mathbf{x} = \twf{f}{0} \wedge \stf{x}{1} ,\quad\quad
f \mathbf{\tilde{x}} = \stf{f}{0} \wedge \twf{x}{1}, \quad\quad
\tilde{f} \mathbf{\tilde{x}} = \twf{f}{0} \wedge \twf{x}{1} , \\
&f \mathbf{x}^\perp = \stf{f}{0} \wedge \tstar \stf{x}{1} = \stf{f}{0} \wedge \twf{x}{n-1}, \quad\quad
\tilde{f} \mathbf{x}^\perp = \twf{f}{0} \wedge \tstar \stf{x}{1} = \twf{f}{0} \wedge \twf{x}{n-1}, \\
\quad\quad
& f \mathbf{\tilde{x}}^\perp  = \stf{f}{0} \wedge \tstar \twf{x}{1} = \stf{f}{0} \wedge \stf{x}{n-1} , \quad\quad
\tilde{f} \mathbf{\tilde{x}}^\perp = \twf{f}{0} \wedge \tstar \twf{x}{1} = \twf{f}{0} \wedge \stf{x}{n-1}  ; 
\end{align}
\item dot product in vector calculus $\leftrightarrow$ Hodge star and wedge product in split exterior calculus:
\begin{align}
& \mathbf{x} \cdot \mathbf{y} = \tstar (\stf{x}{1} \wedge \tstar \stf{y}{1}) = \tstar (\stf{x}{1} \wedge \twf{y}{n-1}) ,  \quad\quad
\mathbf{x} \cdot \mathbf{\tilde{y}} = \tstar (\stf{x}{1} \wedge \tstar \twf{y}{1}) = \tstar (\stf{x}{1} \wedge \stf{y}{n-1}), \quad\quad \\
& \mathbf{\tilde{x}} \cdot \mathbf{y} = \tstar (\twf{x}{1} \wedge \tstar \stf{y}{1}) = \tstar (\twf{x}{1} \wedge \twf{y}{n-1}), \quad\quad
\mathbf{\tilde{x}} \cdot \mathbf{\tilde{y}} = \tstar (\twf{x}{1} \wedge \tstar \twf{y}{1}) = \tstar (\twf{x}{1} \wedge \stf{y}{n-1}) ;
\end{align}
\item inner product in vector calculus $\leftrightarrow$ inner product in split exterior calculus:
\begin{align}
&(f,g) = \left< \stf{f}{0}, \stf{g}{0} \right> = \left< \twf{f}{n}, \twf{g}{n} \right>  ,\quad\quad
(\tilde{f},\tilde{g}) = \left< \twf{f}{0}, \twf{g}{0} \right> = \left< \stf{f}{n}, \stf{g}{n} \right> , \\
&(\mathbf{x}, \mathbf{y}) = \left<\stf{x}{1}, \stf{y}{1}\right> = \left<\twf{x}{n-1}, \twf{y}{n-1}\right> ,  \quad\quad
(\mathbf{\tilde{x}}, \mathbf{\tilde{y}}) = \left<\twf{x}{1}, \twf{y}{1}\right> = \left<\stf{x}{n-1}, \stf{y}{n-1}\right> .
\end{align}
\end{itemize}

Using these identities is one way to derive the split exterior calculus form of equations from vector calculus.
\end{document}